# A novel X-FEM based fast computational method for crack propagation

Zhenxing Cheng[a], Hu Wang[a*]

[a] *State Key Laboratory of Advanced Design and Manufacturing for Vehicle Body, Hunan University, Changsha, 410082, P.R. China*

**Abstract** This study suggests a fast computational method for crack propagation, which is based on the extended finite element method (X-FEM). It is well known that the X-FEM might be the most popular numerical method for crack propagation. However, with the increase of complexity of the given problem, the size of FE model and the number of iterative steps are increased correspondingly. To improve the efficiency of X-FEM, an efficient computational method termed decomposed updating reanalysis (DUR) method is suggested. For most of X-FEM simulation procedures, the change of each iterative step is small and it will only lead a local change of stiffness matrix. Therefore, the DUR method is proposed to predict the modified response by only calculating the changed part of equilibrium equations. Compared with other fast computational methods, the distinctive characteristic of the proposed method is to update the modified stiffness matrix with a local updating strategy, which only the changed part of stiffness matrix needs to be updated. To verify the performance of the DUR method, several typical numerical examples have been analyzed and the results demonstrate that this method is a highly efficient method with high accuracy.

* Corresponding author Tel: +86 0731 88655012; fax: +86 0731 88822051.

E-mail address: wanghu@hnu.edu.cn (Hu Wang)

# 1 Introduction

A great amount of engineering practice indicates that the quality and stability of engineering structures are closely related to the internal crack propagation[1, 2]. Therefore, prediction of the path of crack propagation and analysis of the stability of crack are significant for estimating the safety and the reliability of engineering structures. There are many numerical methods have been developed to simulate crack propagation process, such as finite element method (FEM) [3], boundary element method (BEM) [4], meshless method [5, 6], edge-based finite element method (ES-FEM) [7, 8], numerical manifold method (NMM) [9, 10], extended finite element method (X-FEM) [11-13] and so on[14]. Compared with above methods, the X-FEM might be the most popular numerical method for the crack propagation simulation due to its superiority of modeling both strong and weak discontinuities within a standard FE framework. The X-FEM was first developed by Belytschko and Black [15]. They analyzed the crack propagation problem with minimal re-meshing. Then, Dolbow et al [16]. and Moës et al. [17] improved this method by using a Heaviside function to enrichment function, and it also has been extended to 3D static crack modeling by Sukumar et al.[18]. Sequentially, the X-FEM was significantly improved by coupling with the level set method (LSM) which is used to track both the crack position and tips [19]. Moreover, the X-FEM has been applied to multiple engineering fields, such as dynamic crack propagation or branching [20, 21], crack propagation in composites [22] or shells [23-25], multi-field problems [26], multi-material problems [27, 28], solidification [29], shear bands [30], dislocations [31] and so on. More details of the development of X-FEM can be found in Refs [12, 32-34].

Generally, in order to improve the accuracy of simulation, a very refined mesh with a very small increment of crack propagation or fatigue cycles should be engaged. Correspondingly, the computational cost is expensive. Therefore, reanalysis algorithm is used to improve the efficiency.

Reanalysis, as a fast computational method, is used to predict the response of modified structures efficiently without full analysis, and reanalysis method can be divided into two categories: direct methods (DMs) and approximate methods. DMs can update the inverse of modified stiffness matrix quickly by Sherman-Morrison-Woodbury lemma [35, 36] and obtain the exact response of the modified structure, but usually it can only solve the problems of local or low-rank modifications. In recent decades, many DMs have been suggested. Such as, Song et al. suggested a novel direct reanalysis algorithm based on the binary tree characteristic to update the triangular factorization in sparse matrix solution [37]. Liu et al. applied Cholesky factorization to structural reanalysis [38]. Huang and Wang suggested an independent coefficient (IC) method for large-scale problems with local modification [39]. Compared with DMs, approximate methods can solve the high-rank modifications, but the exact response usually cannot be obtained. The approximate methods mainly include local approximations (LA), global approximations (GA) [40], iterative approximations (IA) [41] and combined approximations (CA) [42-46]. Moreover, many other reanalysis methods have been proposed in recent years. Zuo *et al.* combined reanalysis method with genetic algorithm (GA) [47]. Sun *et al.* extended the reanalysis method into a structural optimization process [48]. To improve the efficiency of reanalysis method, He *et al.* developed a multiple-GPU based parallel IC reanalysis method [49]. Materna et al. applied the reanalysis method to nonlinear problems [50].

In this study, a fast computational method is proposed to model crack propagation under the framework of X-FEM. It is observed that modeling crack propagation by the X-FEM will bring the additional DOFs in each iterative step and it will lead to a local change of stiffness matrix. Considering this characteristic, the decomposed updating reanalysis (DUR) method theoretically can improve the efficiency of X-FEM significantly. Moreover, a local updating stiffness matrix strategy is suggested to improve the efficiency of stiffness matrix assembling. Furthermore, in order to guarantee the accuracy and efficiency of the DUR method, a local strategy has been introduced to update the Cholesky factorization of stiffness matrix efficiently.

The rest of this paper is organized as follows. The basic theories of X-FEM are briefly introduced in Section 2. The details of DUR method are described in Section 3. Then, several numerical examples are tested in Section 4 to investigate the performance of the DUR method. Finally, conclusions are summarized in Section 5.

## 2 Basics theories of XFEM

### 2.1 X-FEM approximation

In the X-FEM, the standard FEM shape function should be enriched by the enrichment function. Assume that the enriched displacement approximation of X-FEM can be defined as:

$$\mathbf{u}^h(\mathbf{x}) = \underbrace{\sum_{I\in\Omega} N_I(\mathbf{x})\mathbf{u}_I}_{\mathbf{u}^{standard}} + \underbrace{\sum_{J\in\Omega_E} \psi(\mathbf{x}) N_J(\mathbf{x})\mathbf{q}_J}_{\mathbf{u}^{enrich}}, \quad (1)$$

where $N_I$ and $\mathbf{u}_I$ denote the standard FEM shape function and nodal degrees of freedom (DOF), respectively. The enriched displacement approximation should be divided into two parts: the standard finite element approximation and partition of unity enrichment. The $\psi(\mathbf{x})$ means enrichment function while the $\mathbf{q}_J$ is the additional nodal degrees of freedom.

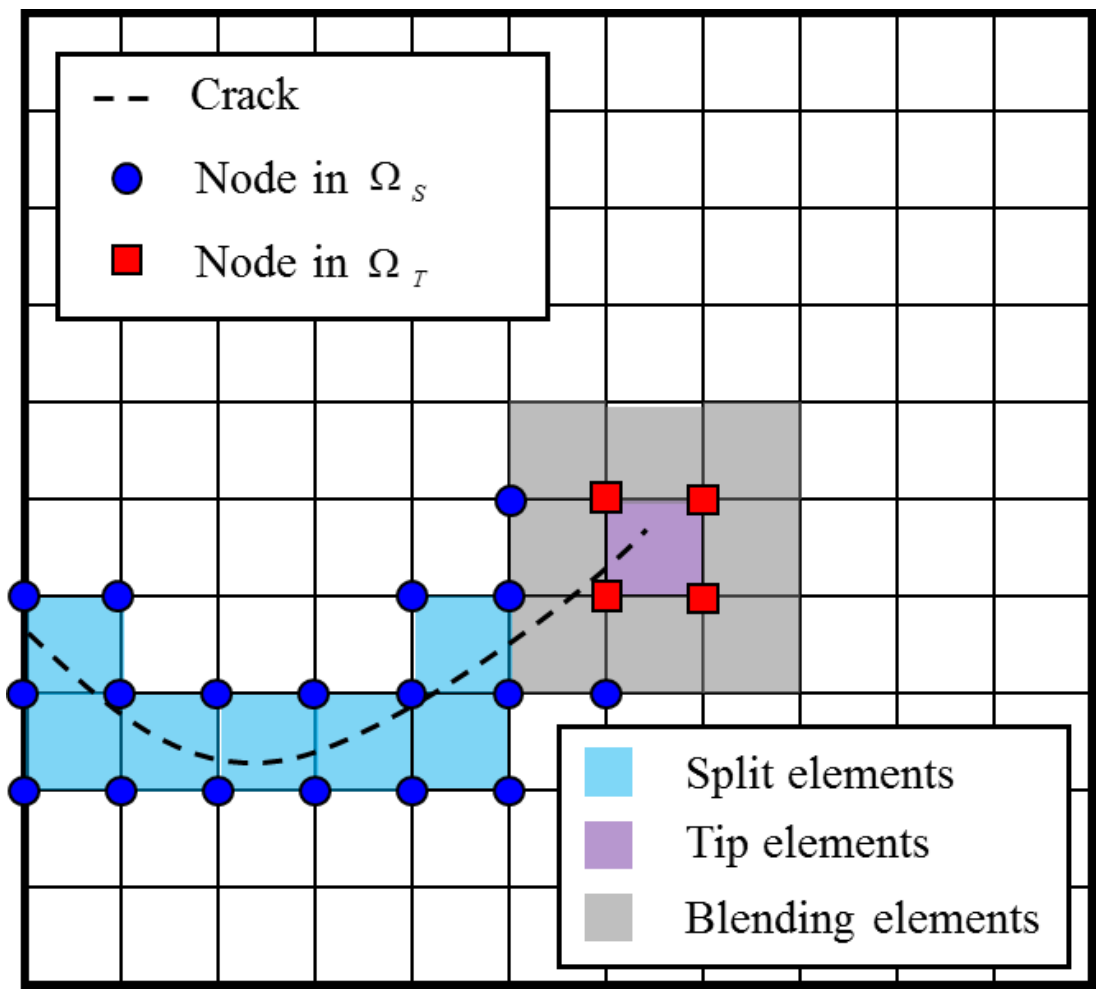

Fig. 1 An arbitrary crack line in a structured mesh

Consider an arbitrary crack in a structured mesh as shown in Fig. 1, then Eq.(1) can be

rewritten as

$$\mathbf{u}^h(\mathbf{x}) = \sum_{I\in\Omega} N_I(\mathbf{x})\mathbf{u}_I + \sum_{I\in\Omega_S} H_I(\mathbf{x})\ N_I(\mathbf{x})\mathbf{a}_I + \sum_{I\in\Omega_T}\sum_{\alpha=1}^{4}\Phi_{I,\alpha}(\mathbf{x})\ N_I(\mathbf{x})\mathbf{b}_I^{\alpha}\,, \tag{2}$$

where $\Omega$ is the solution domain, $\Omega_S$ is the domain cut by crack, $\Omega_T$ is the domain which crack tip located, $H(\mathbf{x})$ is the shifted Heaviside enrichment and $\Phi_\alpha(\mathbf{x})$ is the shifted crack tip enrichment. The details of $H(\mathbf{x})$ and $\Phi_\alpha(\mathbf{x})$ are given as following:

$$H(\mathbf{x}) = \begin{cases} +1 & Above\ \ crack \\ -1 & Below\ \ crack \end{cases}, \tag{3}$$

$$\{\Phi_\alpha(\mathbf{x})\}_{\alpha=1}^{4} = \sqrt{\mathbf{r}}\left\{\sin\frac{\theta}{2}, \cos\frac{\theta}{2}, \sin\theta\sin\frac{\theta}{2}, \sin\theta\cos\frac{\theta}{2}\right\}. \tag{4}$$

The discrete X-FEM equations are obtained by substituting Eq.(2) into the principle of virtual work. Assume that the discrete equations can be defined as

$$\begin{bmatrix} \mathbf{K}_{uu} & \mathbf{K}_{ua} & \mathbf{K}_{ub} \\ \mathbf{K}_{ua}^T & \mathbf{K}_{aa} & \mathbf{K}_{ab} \\ \mathbf{K}_{ub}^T & \mathbf{K}_{ab}^T & \mathbf{K}_{bb} \end{bmatrix} \begin{Bmatrix} \mathbf{u} \\ \boldsymbol{a} \\ \mathbf{b} \end{Bmatrix} = \begin{Bmatrix} \mathbf{F}_u \\ \mathbf{F}_a \\ \mathbf{F}_b \end{Bmatrix}, \tag{5}$$

where $\mathbf{K}_{uu}$ is the traditional finite element stiffness matrix, $\mathbf{K}_{ua}, \mathbf{K}_{aa}, \mathbf{K}_{ab}$ are components with Heaviside enrichment and $\mathbf{K}_{ub}, \mathbf{K}_{ab}, \mathbf{K}_{bb}$ are components with crack tip enrichment.

## 2.2 Crack propagation model

Generally, the direction and magnitude of crack propagation at each iterative step are used to determine how the crack will propagate. The direction of crack propagation is found from the maximum circumferential stress criterion and the crack will propagate in the direction where $\sigma_{\theta\theta}$ is maximum [51]. The angle of crack propagation is defined as

$$\theta = 2\arctan\frac{1}{4}\left(\frac{K_I}{K_{II}} - sign\ K_{II}\sqrt{\left(\frac{K_I}{K_{II}}\right)^2 + 8}\right), \tag{6}$$

where $\theta$ is defined in the crack tip coordinate system, $K_I$ and $K_{II}$ are the mixed-mode stress intensity factors. The details are given in the reference [51].

There are two main quasi-static manners when modeling crack growth. The first one assumes a constant increment of crack growth at each cycle [16]while the other option is to assume a constant number of cycles and apply a fatigue crack growth law to predict the crack growth increment for the fixed number of cycles [52]. In this study, a fixed increment of crack growth $\Delta a$ is considered.

# 3 Decomposed updating reanalysis method

## 3.1 Framework of the DUR method

The DUR method is proposed to model quasi-static crack propagation under the framework of X-FEM and the framework of the DUR method is presented in Fig. 2.

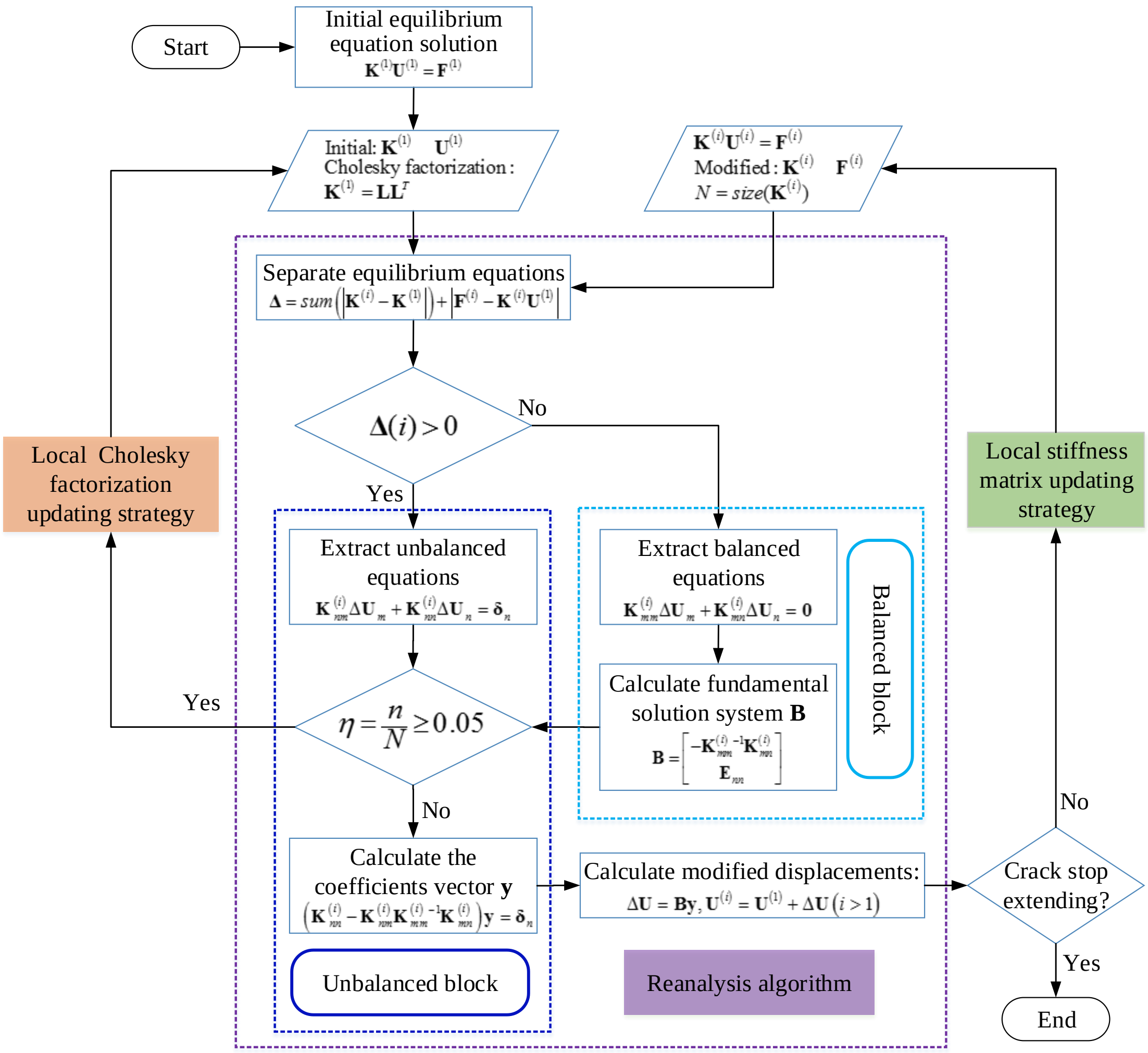

Fig. 2 The framework of the DUR method

It can be found that the DUR method includes three parts: local updating stiffness matrix strategy, local Cholesky factorization updating strategy and reanalysis algorithm. The local stiffness matrix updating strategy is suggested to improve the efficiency of stiffness matrix assembling according to the characteristic of X-FEM. Moreover, considering the local change of stiffness matrix, the specific reanalysis algorithm is used to predict the modified response and improve the efficiency of the X-FEM. Furthermore, in order to guarantee the accuracy and the efficiency of the suggested reanalysis method, a local strategy has been introduced to update the Cholesky factorization of stiffness matrix efficiently. More details can be found in Section 3.2, 3.3 and 3.4.

### 3.2 Local stiffness matrix updating strategy

As mentioned above, the stiffness matrix can be given as the following form, and each component has been marked by different colors associated with Fig. 1.

$$\mathbf{K} = \begin{bmatrix} \mathbf{K}_{uu} & \mathbf{K}_{ua} & \mathbf{K}_{ub} \\ \mathbf{K}_{ua}^{T} & \mathbf{K}_{aa} & \mathbf{K}_{ab} \\ \mathbf{K}_{ub}^{T} & \mathbf{K}_{ab}^{T} & \mathbf{K}_{bb} \end{bmatrix}$$

Fig. 3 Color-marked stiffness matrix

Consider a small increment of crack based on Fig. 1, and the comparison is shown in Fig. 4.

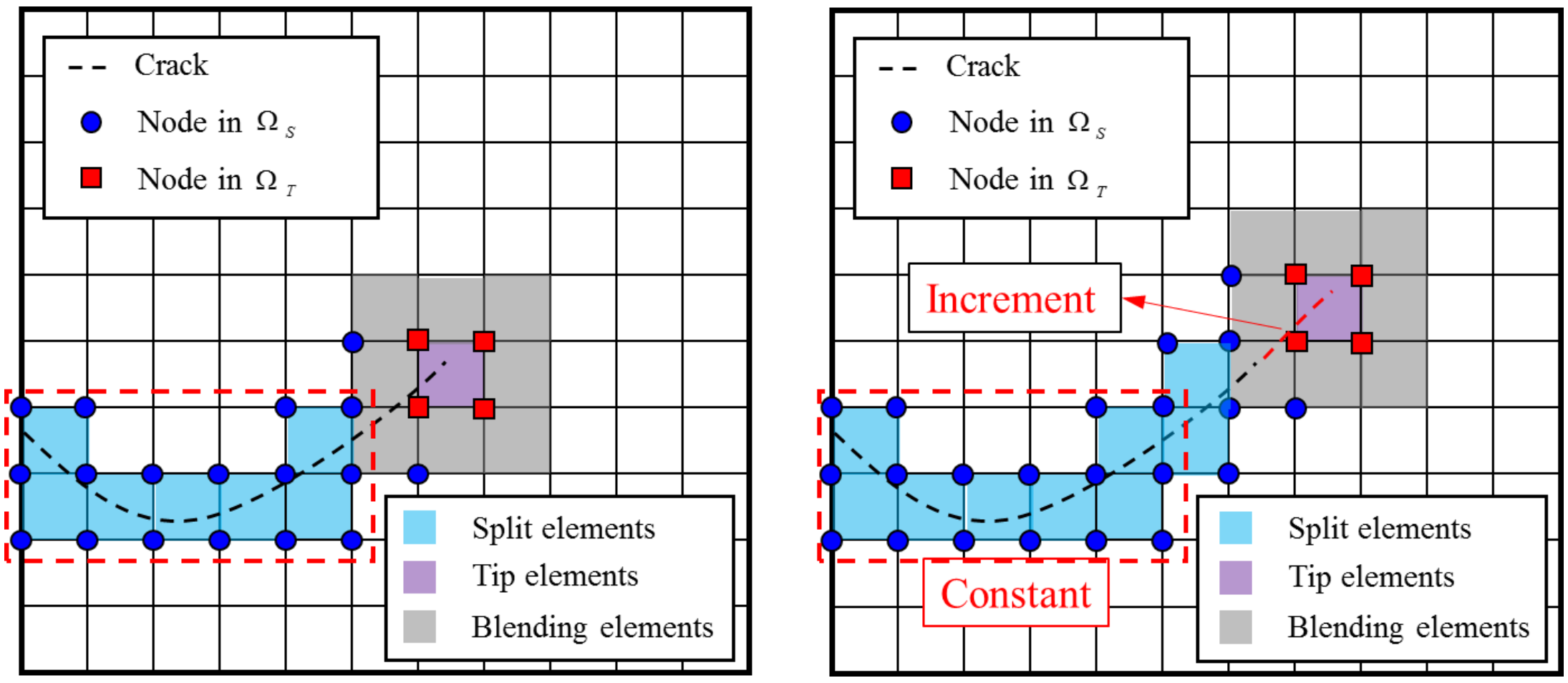

Fig. 4 A small increment of crack propagation

It can be found that the stiffness component $\mathbf{K}_{uu}$ associated with the standard finite element approximation will be constant at each iterative step of crack propagation. It means that the change part of stiffness matrix is the enriched part, which should be much smaller than the un-enriched part. Furthermore, it can be found that once an element had become a split element, it will be a split element during entire computational procedure. It implies that once an element has been enriched by Heaviside function, the value of stiffness matrix will be a constant in the subsequent iterative steps. Therefore, the changed part of the stiffness matrix in each iterative step will be a small part as shown in Fig. 5.

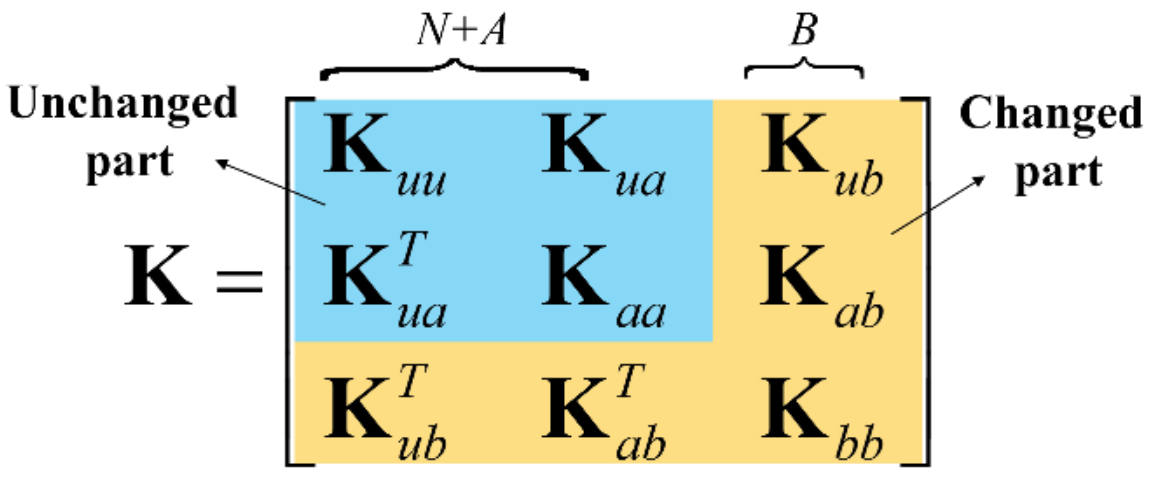

Fig. 5 Changed and unchanged part of stiffness matrix

Assume that the unchanged part is a $(N+A)\times(N+A)$ matrix and the changed part contains a $B\times B$ matrix with $(N+A)\times B$, $B\times(N+A)$ coupling matrices. If $N>>A+B$, only a small part of stiffness matrix needs to be updated in each iterative steps. Obviously, it should be high efficient and the matrix assembling cost should be significantly reduced.

### 3.3 The reanalysis algorithm for X-FEM

The reanalysis algorithm is used to predict the response of the current iterative step by using the information of the first iterative step. In this study, a reanalysis algorithm is proposed for X-FEM according to the characteristic of X-FEM. The reanalysis algorithm avoids the full analysis after the first iterative step, and the response of the subsequent iterative steps can be efficiently obtained. The details of the reanalysis algorithm are described as following.

Assume that the equilibrium equation of the *i-th* iterative step is

$$\mathbf{K}^{(i)}\mathbf{U}^{(i)}=\mathbf{F}^{(i)}, \tag{7}$$

where $\mathbf{U}^{(i)}$ is the displacement in the *i-th* iterative step, and the equilibrium equation of the first iterative step can be given as

$$\mathbf{K}^{(1)}\mathbf{U}^{(1)}=\mathbf{F}^{(1)}. \tag{8}$$

Assume that the solution of the equation in the *i-th* iterative step can be defined as

$$\mathbf{U}^{(i)}=\mathbf{U}^{(1)}+\Delta\mathbf{U}\left(i>1\right), \tag{9}$$

then substitute Eq.(9) into Eq.(7),

$$\mathbf{K}^{(i)}\left(\mathbf{U}^{(1)}+\Delta\mathbf{U}\right)=\mathbf{F}^{(i)}\left(i>1\right), \tag{10}$$

or

$$\mathbf{K}^{(i)}\Delta\mathbf{U}=\mathbf{F}^{(i)}-\mathbf{K}^{(i)}\mathbf{U}^{(1)}\left(i>1\right). \tag{11}$$

Define the residual value of displacement $\boldsymbol{\delta}$ as

$$\boldsymbol{\delta}=\mathbf{F}^{(i)}-\mathbf{K}^{(i)}\mathbf{U}^{(1)}\left(i>1\right), \tag{12}$$

then Eq.(11) can be written as

$$\mathbf{K}^{(i)}\Delta\mathbf{U}=\boldsymbol{\delta}. \tag{13}$$

Consider that only a small part of stiffness matrix will change in every iterative step, the most part of $\boldsymbol{\delta}$ should be zero. Based on this property, the $\mathbf{U}^{(i)}$ can be divided into two blocks: unbalanced and balanced blocks, according to Eq.(14):

$$\boldsymbol{\Delta}=sum\left(\left|\mathbf{K}^{(i)}-\mathbf{K}^{(1)}\right|\right)+\left|\boldsymbol{\delta}\right|. \tag{14}$$

If $\left|\boldsymbol{\Delta}\left(j\right)\right|>0$, the *j-th* DOF is unbalanced, otherwise the *j-th* DOF is balanced. Accord to this, the Eq.(13) can be rewritten as

$$\begin{bmatrix}\mathbf{K}_{mm}^{(i)} & \mathbf{K}_{mn}^{(i)}\\ \mathbf{K}_{nm}^{(i)} & \mathbf{K}_{nn}^{(i)}\end{bmatrix}\begin{Bmatrix}\Delta\mathbf{U}_m\\ \Delta\mathbf{U}_n\end{Bmatrix}=\begin{Bmatrix}\mathbf{0}\\ \boldsymbol{\delta}_n\end{Bmatrix}, \tag{15}$$

where *m* is the number of balanced DOFs, and *n* is the number of unbalanced DOFs. Equation (15) can be rewritten as

$$\mathbf{K}_{mm}^{(i)}\Delta\mathbf{U}_m+\mathbf{K}_{mn}^{(i)}\Delta\mathbf{U}_n=\mathbf{0} \tag{16}$$

and

$$\mathbf{K}_{nm}^{(i)}\Delta\mathbf{U}_m+\mathbf{K}_{nn}^{(i)}\Delta\mathbf{U}_n=\boldsymbol{\delta}_n. \tag{17}$$

Obviously, Eq.(16) is a homogeneous equation set which has infinite solutions, and the solution $\mathbf{U}^{(1)}$ is one of them.

Let $\Delta\mathbf{U}_n=\left\{0 \ \cdots \ 0 \ 1 \ 0 \ \cdots \ 0\right\}^T$ (where 1 is the *k-th* member of $\Delta\mathbf{U}_n$, and $k=1,2,\ldots,n$), then the fundamental solution system of Eq. (16) can be obtained as

$$\mathbf{B} = \begin{bmatrix} -\mathbf{K}_{mm}^{(i)\ -1}\mathbf{K}_{mn}^{(i)} \\ \mathbf{E}_{nn} \end{bmatrix}, \tag{18}$$

where $\mathbf{E}_{nn}$ is a rank-*n* unit matrix.

Define the general solution of Eq.(16) is

$$\Delta\mathbf{U} = \mathbf{B}\mathbf{y}, \tag{19}$$

where **y** is a dimension-*n* vector.

Then, substitute Eq.(19) into Eq.(17) to find a unique solution of Eq.(13), obtain

$$\left(\mathbf{K}_{nn}^{(i)} - \mathbf{K}_{nm}^{(i)}\mathbf{K}_{mm}^{(i)\ -1}\mathbf{K}_{mn}^{(i)}\right)\mathbf{y} = \boldsymbol{\delta}_n. \tag{20}$$

Solve Eq.(20), **y** can be obtained, and then $\Delta\mathbf{U}$ can be obtained by Eq.(19). Sequentially, the $\mathbf{U}^{(i)}$ can be obtained by Eq.(9). Because only a small part of stiffness matrix will change in each iterative step, the number *n* should not be too large. Therefore, Eq.(20) is a small scale problem. The key point is how to obtain $\mathbf{K}_{mm}^{(i)}$ from the $\mathbf{K}^{(1)}$.

Assume that

$$\mathbf{K}^{(1)} = \begin{bmatrix} \mathbf{K}_{mm}^{(1)} & \mathbf{K}_{mn}^{(1)} \\ \mathbf{K}_{nm}^{(1)} & \mathbf{K}_{nn}^{(1)} \end{bmatrix}, \tag{21}$$

and then the Cholesky factorization of $\mathbf{K}^{(1)}$ can be defined as

$$\begin{aligned} \mathbf{K}^{(1)} = \mathbf{L}\mathbf{L}^T &= \begin{bmatrix} \mathbf{L}_{mm} & \mathbf{0} \\ \mathbf{L}_{nm} & \mathbf{L}_{nn} \end{bmatrix} \begin{bmatrix} \mathbf{L}_{mm}^T & \mathbf{L}_{nm}^T \\ \mathbf{0} & \mathbf{L}_{nn}^T \end{bmatrix} \\ &= \begin{bmatrix} \mathbf{L}_{mm}\mathbf{L}_{mm}^T & \mathbf{L}_{mm}\mathbf{L}_{nm}^T \\ \mathbf{L}_{nm}\mathbf{L}_{mm}^T & \mathbf{L}_{nm}\mathbf{L}_{nm}^T + \mathbf{L}_{nn}\mathbf{L}_{nn}^T \end{bmatrix} \end{aligned}, \tag{22}$$

where $\mathbf{L}_{mm}$ and $\mathbf{L}_{nn}$ are lower triangular matrices.

Compared with Eq.(21), it can be found that

$$\mathbf{K}_{mm}^{(i)} = \mathbf{K}_{mm}^{(1)} = \mathbf{L}_{mm}\mathbf{L}_{mm}^T. \tag{23}$$

Therefore, the fundamental solution system **B** can be calculated by Eq.(18), and the Cholesky factorization of stiffness matrix in the first iterative step can directly re-used. In describe this method clearly, the $\mathbf{K}_{mm}^{(i)}$ can be associated with Fig. 4. Assume that the left figure of Fig. 4 shows the crack position in the first iterative step, and the right figure of Fig. 4 shows the crack position in the second iterative step. Then the stiffness matrices of X-FEM in the first and the second iterative steps can be written as Fig. 6 and Fig. 7 respectively.

$$\mathbf{K}^{(1)} = \begin{bmatrix} \mathbf{K}_{uu}^{(1)} & \mathbf{K}_{ua}^{(1)} & \mathbf{K}_{ub}^{(1)} \\ \mathbf{K}_{ua}^{(1)T} & \mathbf{K}_{aa}^{(1)} & \mathbf{K}_{ab}^{(1)} \\ \mathbf{K}_{ub}^{(1)T} & \mathbf{K}_{ab}^{(1)T} & \mathbf{K}_{bb}^{(1)} \end{bmatrix}$$

Fig. 6 The stiffness matrix of the first iterative step

$$\mathbf{K}^{(2)} = \begin{bmatrix} \mathbf{K}_{uu}^{(1)} & \mathbf{K}_{ua}^{(1)} & \mathbf{K}_{ua}^{(2)} & \mathbf{K}_{ub}^{(2)} \\ \mathbf{K}_{ua}^{(1)T} & \mathbf{K}_{aa}^{(1)} & \mathbf{0} & \mathbf{0} \\ \mathbf{K}_{ua}^{(2)T} & \mathbf{0} & \mathbf{K}_{aa}^{(2)} & \mathbf{K}_{ab}^{(2)} \\ \mathbf{K}_{ub}^{(2)T} & \mathbf{0} & \mathbf{K}_{ab}^{(2)T} & \mathbf{K}_{bb}^{(2)} \end{bmatrix}$$

Fig. 7 The stiffness matrix of the second iterative step

It can be found that the red-marked part is constant as shown in Fig. 6 and Fig. 7. Compared with Eq.(15), the $\mathbf{K}_{mm}^{(2)}$ can be give as

$$\mathbf{K}_{mm}^{(2)} = \begin{bmatrix} \mathbf{K}_{uu}^{(1)} & \mathbf{K}_{ua}^{(1)} \\ \mathbf{K}_{ua}^{(1)T} & \mathbf{K}_{aa}^{(1)} \end{bmatrix}, \tag{24}$$

and the $\mathbf{K}_{mn}^{(2)}$, $\mathbf{K}_{nm}^{(2)}$ and $\mathbf{K}_{nn}^{(2)}$ can be give as

$$\mathbf{K}_{mn}^{(2)} = \begin{bmatrix} \mathbf{K}_{ua}^{(2)} & \mathbf{K}_{ub}^{(2)} \\ \mathbf{0} & \mathbf{0} \end{bmatrix}, \tag{25}$$

$$\mathbf{K}_{nm}^{(2)} = \begin{bmatrix} \mathbf{K}_{ua}^{(2)T} & \mathbf{0} \\ \mathbf{K}_{ub}^{(2)T} & \mathbf{0} \end{bmatrix}, \tag{26}$$

$$\mathbf{K}_{nn}^{(2)} = \begin{bmatrix} \mathbf{K}_{aa}^{(2)} & \mathbf{K}_{ab}^{(2)} \\ \mathbf{K}_{ab}^{(2)T} & \mathbf{K}_{bb}^{(2)} \end{bmatrix}. \tag{27}$$

### 3.4 Local Cholesky factorization updating strategy

Generally, the accuracy and the efficiency of the reanalysis methods are depended on the percentage of changed part, and if the percentage of changed part is too large, the accuracy and the efficiency should be unavailable for engineering problems. Therefore, a suitable critical value of changed percentage should be stipulated, and usually 5% is chosen for the critical value [53].

In the DUR method, the changed part is composed by unbalanced DOFs. Then the changed percentage $\eta$ can be defined as the following form:

$$\eta = \frac{n}{N} \times 100\% , \tag{28}$$

where $n$ is the number of unbalanced DOFs and $N$ is the number of total DOFs. In order to guarantee the accuracy and efficiency of the reanalysis method, a suitable critical value of $\eta$ should be set. In this study, 5% is also chose as the critical value of $\eta$, so the initial information needs to be updating when the $\eta > 5\%$ and the key issue is how to obtain the Cholesky factorization of the initial stiffness matrix efficiently. Therefore, a local updating Cholesky factorization strategy has been suggested. This strategy is based on the property that only a small part of stiffness matrix will change in each iterative step, so that the Cholesky factorization of modified stiffness matrix can be updating based on the Cholesky factorization of initial stiffness matrix. An example is given to explain this strategy, which is based on Fig. 6 and Fig. 7.

Assume that Fig. 6 means the initial stiffness matrix and Fig. 7 means the modified stiffness matrix. Recall the Cholesky factorization of initial stiffness matrix takes the form

$$\mathbf{K}^{(1)} = \begin{bmatrix} \mathbf{K}_{11}^{(1)} & \mathbf{K}_{12}^{(1)} \\ \mathbf{K}_{12}^{(1)T} & \mathbf{K}_{22}^{(1)} \end{bmatrix} = \begin{bmatrix} \mathbf{L}_{11}^{(1)} & \mathbf{0} \\ \mathbf{L}_{21}^{(1)} & \mathbf{L}_{22}^{(1)} \end{bmatrix} \begin{bmatrix} \mathbf{L}_{11}^{(1)T} & \mathbf{L}_{21}^{(1)T} \\ \mathbf{0} & \mathbf{L}_{22}^{(1)T} \end{bmatrix} = \mathbf{L}^{(1)}\mathbf{L}^{(1)T} , \tag{29}$$

where

$$\begin{aligned} \mathbf{L}_{11}^{(1)} &= chol\,\mathbf{K}_{11}^{(1)} \\ \mathbf{L}_{21}^{(1)} &= \mathbf{L}_{11}^{(1)} / \mathbf{K}_{12}^{(1)} \\ \mathbf{L}_{22}^{(1)} &= chol\,\mathbf{K}_{22}^{(1)} - \mathbf{L}_{21}^{(1)}\mathbf{L}_{21}^{(1)T} \end{aligned} \quad . \tag{30}$$

Compared with Fig. 6, $\mathbf{K}_{11}^{(1)}$ can be considered equivalent to

$$\mathbf{K}_{11}^{(1)} = \begin{bmatrix} \mathbf{K}_{uu}^{(1)} & \mathbf{K}_{ua}^{(1)} \\ \mathbf{K}_{ua}^{(1)T} & \mathbf{K}_{aa}^{(1)} \end{bmatrix}, \tag{31}$$

and $\mathbf{K}_{12}^{(1)}$, $\mathbf{K}_{22}^{(1)}$ can be written as

$$\mathbf{K}_{12}^{(1)} = \begin{bmatrix} \mathbf{K}_{ub}^{(1)} \\ \mathbf{K}_{ab}^{(1)} \end{bmatrix}, \ \mathbf{K}_{22}^{(1)} = \mathbf{K}_{bb}^{(1)}. \tag{32}$$

Then the Cholesky factorization of modified stiffness matrix $\mathbf{L}^{(2)}$ can be obtained by

$$\mathbf{L}^{(2)} = \begin{bmatrix} \mathbf{L}_{11}^{(2)} & \mathbf{0} \\ \mathbf{L}_{21}^{(2)} & \mathbf{L}_{22}^{(2)} \end{bmatrix}, \tag{33}$$

where

$$\begin{aligned} \mathbf{L}_{11}^{(2)} &= chol\,\mathbf{K}_{11}^{(2)} = chol\,\mathbf{K}_{11}^{(1)} = \mathbf{L}_{11}^{(1)} \\ \mathbf{L}_{21}^{(2)} &= \mathbf{L}_{11}^{(2)} / \mathbf{K}_{12}^{(2)} \\ \mathbf{L}_{22}^{(2)} &= chol\,\mathbf{K}_{22}^{(2)} - \mathbf{L}_{21}^{(2)}\mathbf{L}_{21}^{(2)T} \end{aligned} \quad , \tag{34}$$

and compared with Fig. 7,

$$\mathbf{K}_{12}^{(2)} = \begin{bmatrix} \mathbf{K}_{ua}^{(2)} & \mathbf{K}_{ub}^{(2)} \\ \mathbf{0} & \mathbf{0} \end{bmatrix}, \ \mathbf{K}_{22}^{(2)} = \begin{bmatrix} \mathbf{K}_{aa}^{(2)} & \mathbf{K}_{ab}^{(2)} \\ \mathbf{K}_{ab}^{(2)T} & \mathbf{K}_{bb}^{(2)} \end{bmatrix}. \tag{35}$$

Then $\mathbf{L}_{11}^{(1)}$ can be reused at subsequent iterative steps of crack propagation, and usually the size of $\mathbf{L}_{11}^{(1)}$ is much large than other parts. Therefore, the strategy should save much computational cost than calculate the Cholesky factorization of the entire stiffness matrix directly.

## 4 Numerical examples

In order to test the accuracy and efficiency of the DUR method, three examples are tested by the proposed methods. These three cases involve edge and center crack propagation, concentrated and uniformed load problems, thus the performance of the DUR method could be verified thoroughly. In this study, the comparison has been made between the DUR and full analysis, and the errors of displacement, Von Mises stress and Von Mises strain are defined by the following formulas:

$$E_u = \frac{\left\| \mathbf{U}_{PFR} - \mathbf{U}_{FA} \right\|}{\left\| \mathbf{U}_{FA} \right\|} \times 100\% \tag{36}$$

$$E_\sigma = \frac{\left\| \boldsymbol{\sigma}_{PFR} - \boldsymbol{\sigma}_{FA} \right\|}{\left\| \boldsymbol{\sigma}_{FA} \right\|} \times 100\% \tag{37}$$

where $\mathbf{U}_{PFR}$, $\boldsymbol{\sigma}_{PFR}$ mean the results of DUR method, and $\mathbf{U}_{FA}$, $\boldsymbol{\sigma}_{FA}$ mean the results of full analysis.

Moreover, in order to investigate the performance of the DUR method, the CPU running time has been recorded and all the simulations were performed on an Intel(R) Core(TM) i7-5820K 3.30GHz CPU with 32GB of memory within MATLAB R2016b in x64 Windows 7.

### 4.1 Edge crack in a plate with a hole

The problem shown in Fig. 8 is an adaptation of an example presented in reference [54]. The initial crack length is $a_0 = 10mm$, the force $F = 2 \times 10^4 N$, and linear elastic material behavior is assumed. The material is aluminum 7075-T6 with $E = 7.17 \times 10^4 MPa$, $\nu = 0.33$ and a plane strain state is considered. The increment of propagation $\Delta a = 1mm$. For this problem, the result of experimental test with specimens 16 *mm* thick was given by Giner et al. [54] as shown in Fig. 9.

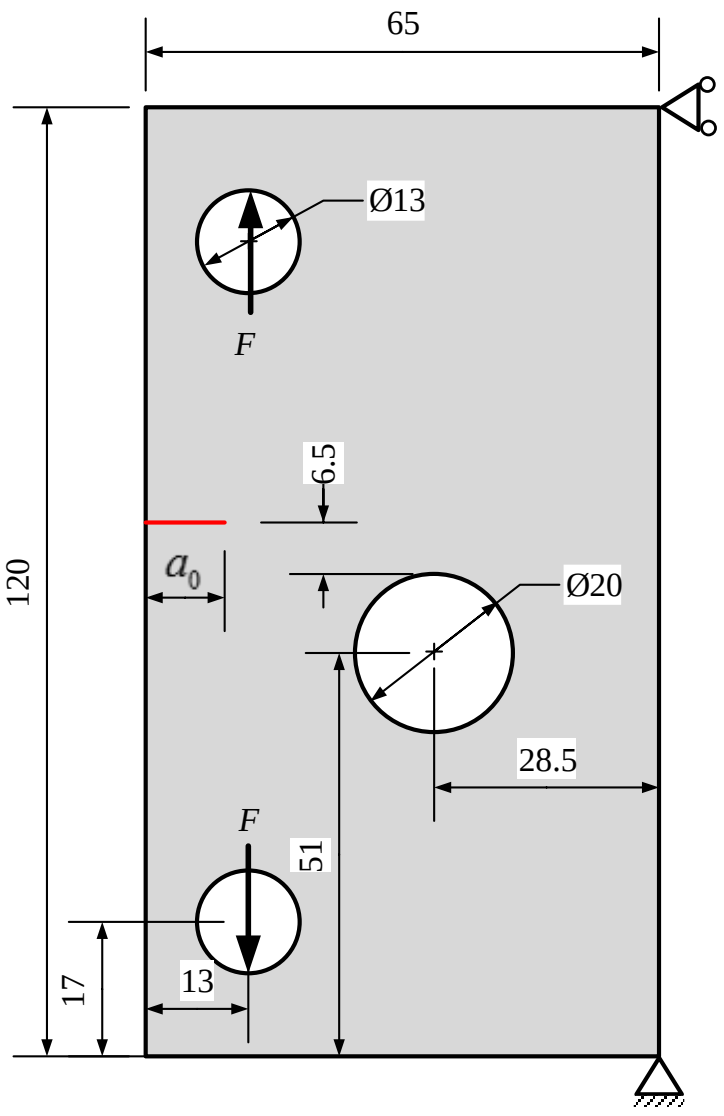

Fig. 8 The geometry of edge crack in a plate with a hole

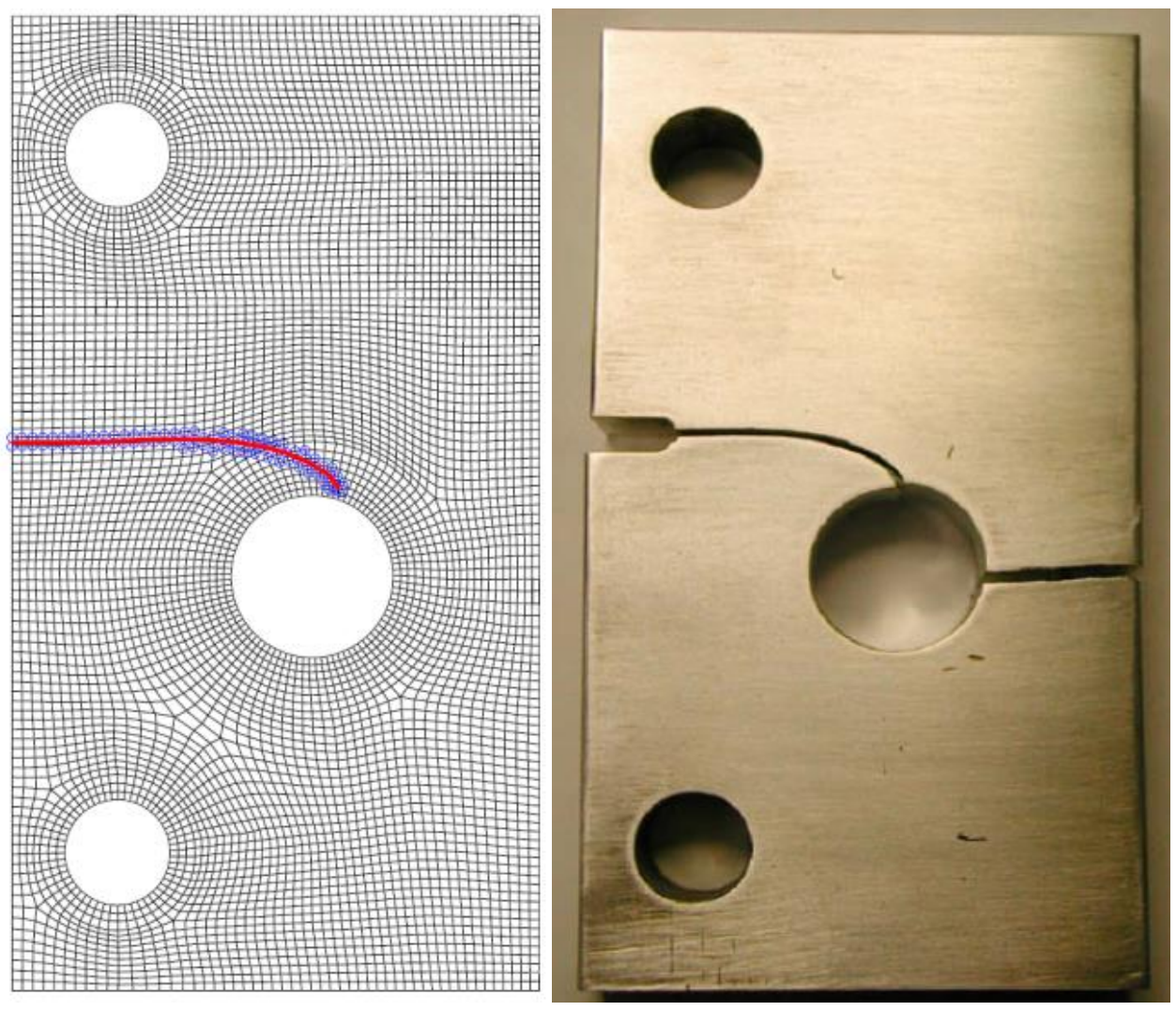

Fig. 9 The result of edge crack in a plate with a hole in the reference [54]

Then solve this case by the DUR method and only the first iterative step need to be calculated by full analysis method while other iterative steps should be predicted efficiently by the DUR method. In order to investigate the accuracy of the DUR method, The comparisons of displacement and stress between the DUR and full analysis are illustrated in Fig. 10 and Fig. 11, respectively. It is obviously that the results of DUR method and full analysis are almost the same. Compared with Fig. 9, it is proved that the DUR method is accurate. Moreover, the errors of each iterative step which is defined by Eq.(36) and Eq.(37) are shown in Fig. 12. It can be found that the DUR method is accurate. The computational costs of the DUR and full analysis are also listed in Tab.

1. It shows that the computational cost of the DUR method is much cheaper than the full analysis method. Moreover, the computational results of crack tip coordinates are shown in Tab. 2 and the results of the DUR and full analysis are also the same.

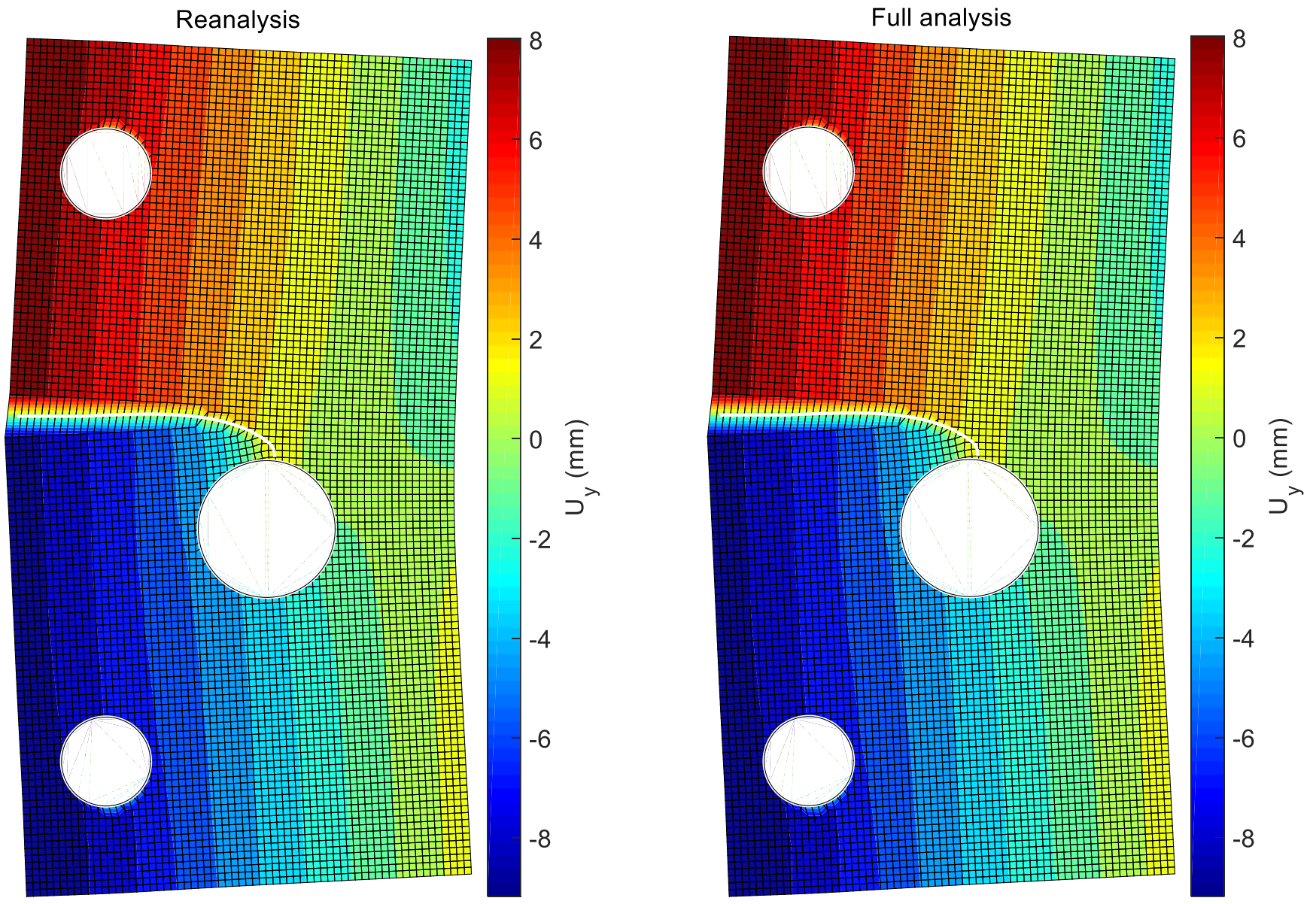

Fig. 10 Displacement contour of edge crack in a plate with a hole

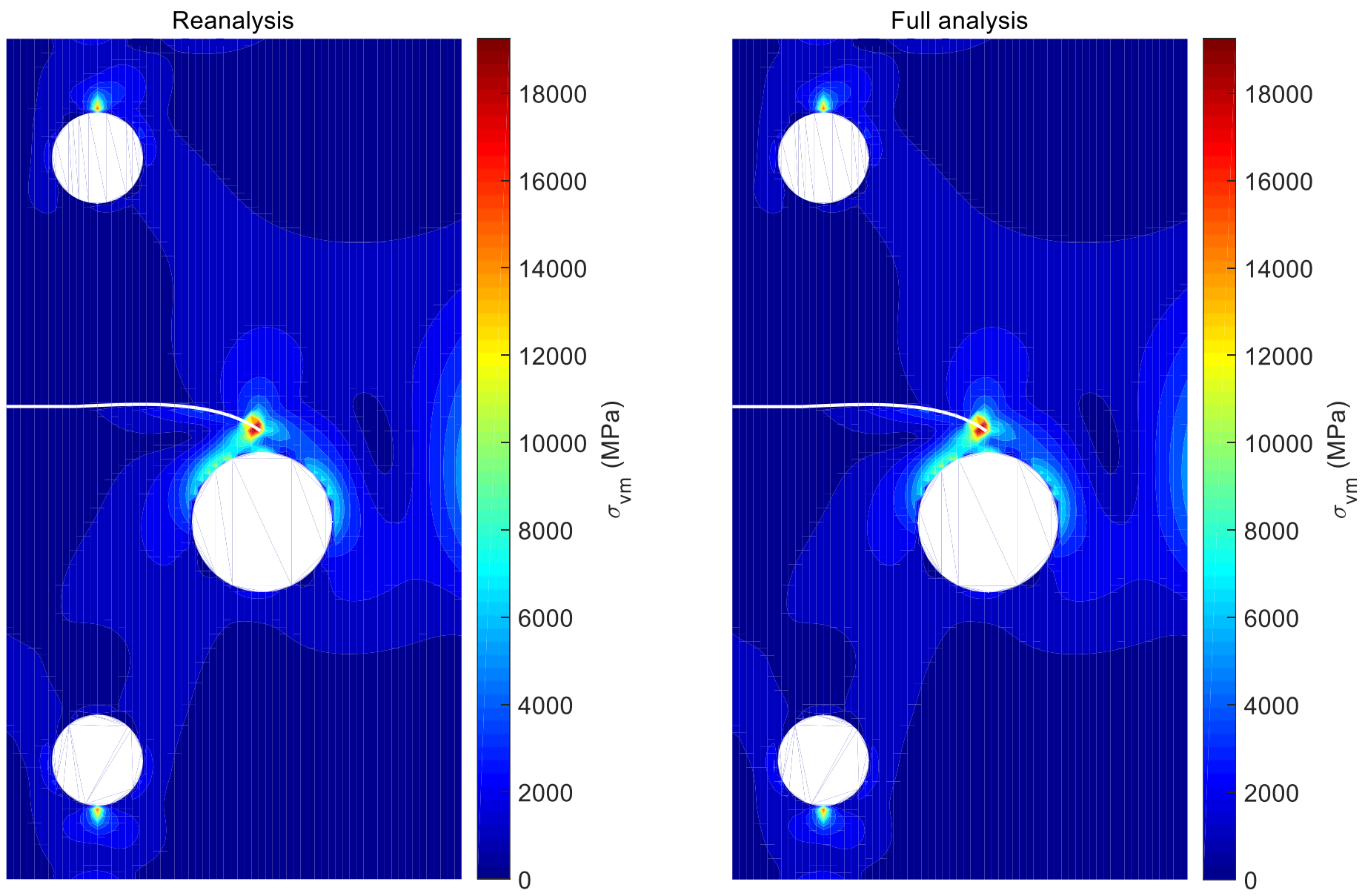

Fig. 11 Von Mises stress contour of edge crack in a plate with a hole

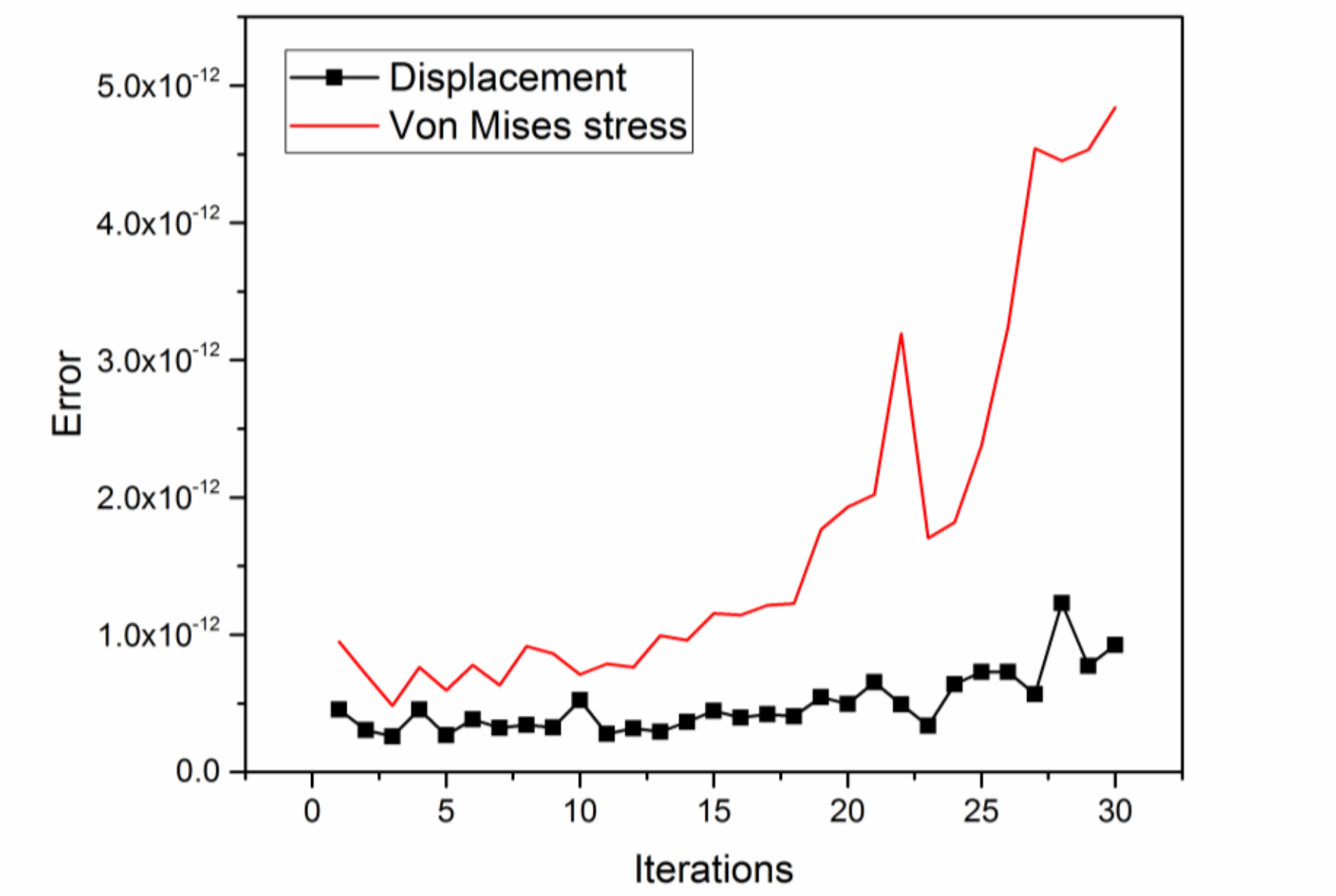

Fig. 12 The errors of each iterative step

Tab. 1 Performance comparison of edge crack in a plate with a hole

| Computational Time/s | | Average Errors | |
|---|---|---|---|
| DUR | Full analysis | Displacement | Von Mises Stress |
| 0.968 | 18.953 | 4.8916e-13 | 1.7356e-12 |

Tab. 2 Crack tip coordinates of edge crack in a plate with a hole

| Iterative steps | Crack tip coordinates of DUR Method | | Crack tip coordinates of full analysis | |
|---|---|---|---|---|
| | X-coordinate | Y-coordinate | X-coordinate | Y-coordinate |
| 0 | 0.000 | 67.500 | 0.000 | 67.500 |
| | 10.000 | 67.500 | 10.000 | 67.500 |
| 1 | 10.998 | 67.559 | 10.998 | 67.559 |
| 5 | 14.994 | 67.727 | 14.994 | 67.727 |
| 10 | 19.993 | 67.805 | 19.993 | 67.805 |
| 15 | 24.986 | 67.560 | 24.986 | 67.560 |
| 20 | 29.918 | 66.766 | 29.918 | 66.766 |
| 25 | 34.605 | 65.059 | 34.605 | 65.059 |
| 30 | 37.695 | 61.523 | 37.695 | 61.523 |

### 4.2 Edge crack in a plate with a circular inclusion

As shown in Fig. 13, a plate with a circular inclusion is considered. Assume that there is an edge crack in the left edge of the plate. The aluminum 7075-T6 is used as the material of plate in the previous case, the Young's modulus $E_1 = 7.17 \times 10^4 MPa$, the Poisson's ratio $\nu = 0.33$, and a plane strain state is assuming. The material of inclusion is considered as carbon fiber reinforced composite and assume that the Young's modulus $E_2 = 2.1\times10^5 MPa$, the Poisson's is the same as the plate. The initial crack length $a_0 = 10mm$ , and the uniformed load $q = 50N/mm$ . The increment of propagation $\Delta a = 1mm$. Then the DUR method is used to solve this problem and only the first iterative step needs to be calculated by full analysis method while others iterative steps should be predicted efficiently by the DUR method.

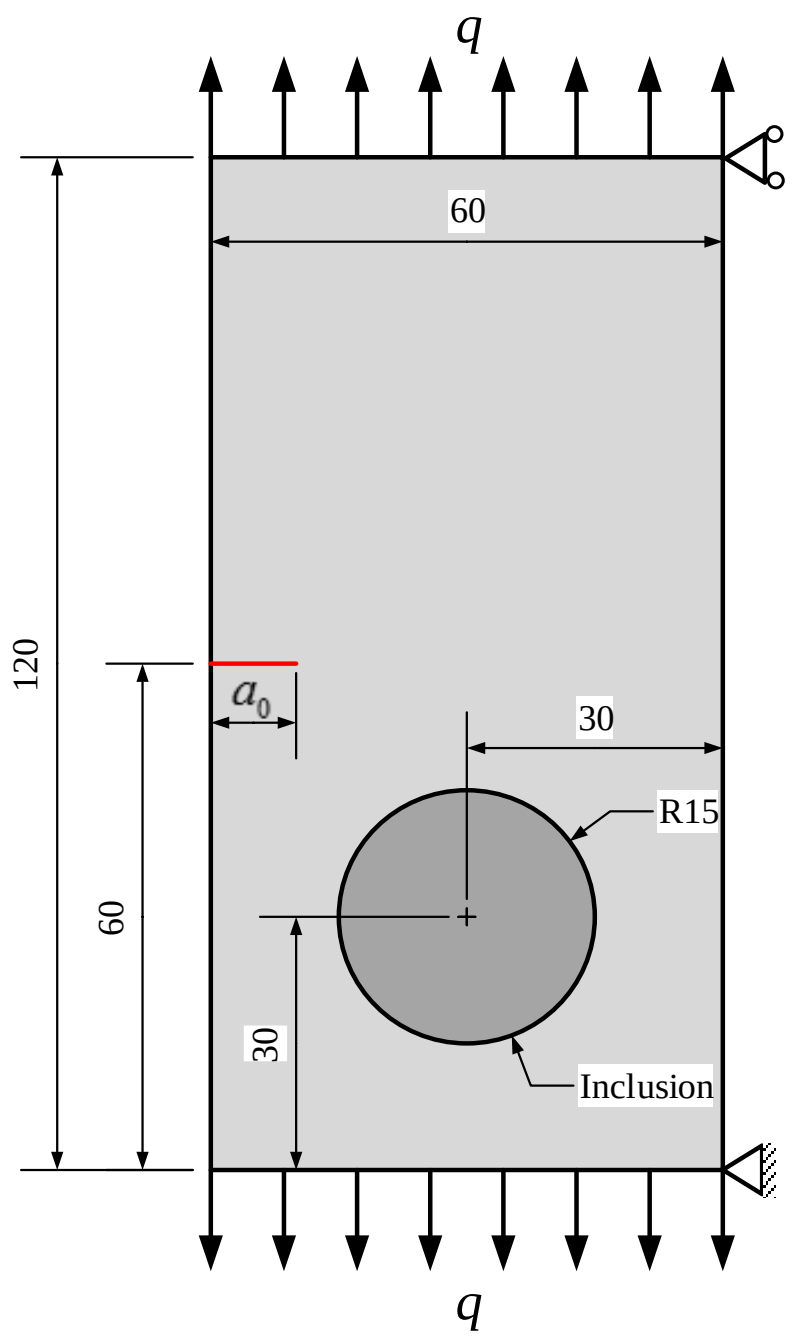

Fig. 13 The geometry of edge crack in a plate with a circular inclusion

The comparisons of displacement and stress are shown as Fig. 14 and Fig. 15, Fig. 15 presents the Von Mises stress results of iterative step 1, 20 and 40. It is obviously that the result of DUR method is very close to the full analysis method. Furthermore, an error list of some selected iterative steps are shown as Fig. 16. It can be found that the DUR method is highly accurate because the maximum of error is $1.0\times10^{-10}$. Moreover, the computational costs of the DUR and full analysis are also listed in Tab. 3. It shows that the efficiency of the DUR method is much higher than the full analysis method. Moreover, the crack tip coordinates are also listed in Tab. 4 and the results of the DUR and full analysis are matched.

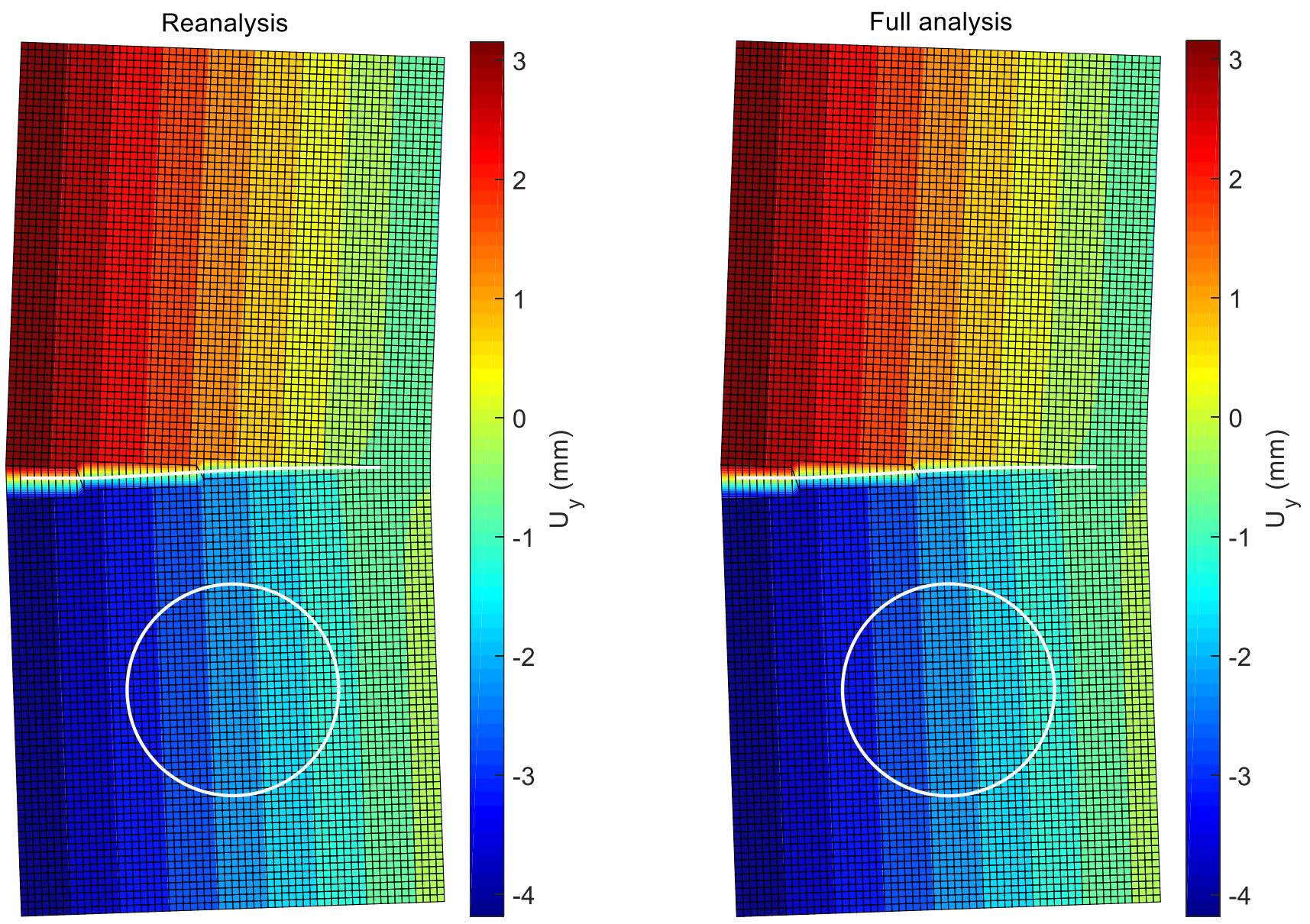

Fig. 14 Displacement contour of edge crack in a plate with a circular inclusion

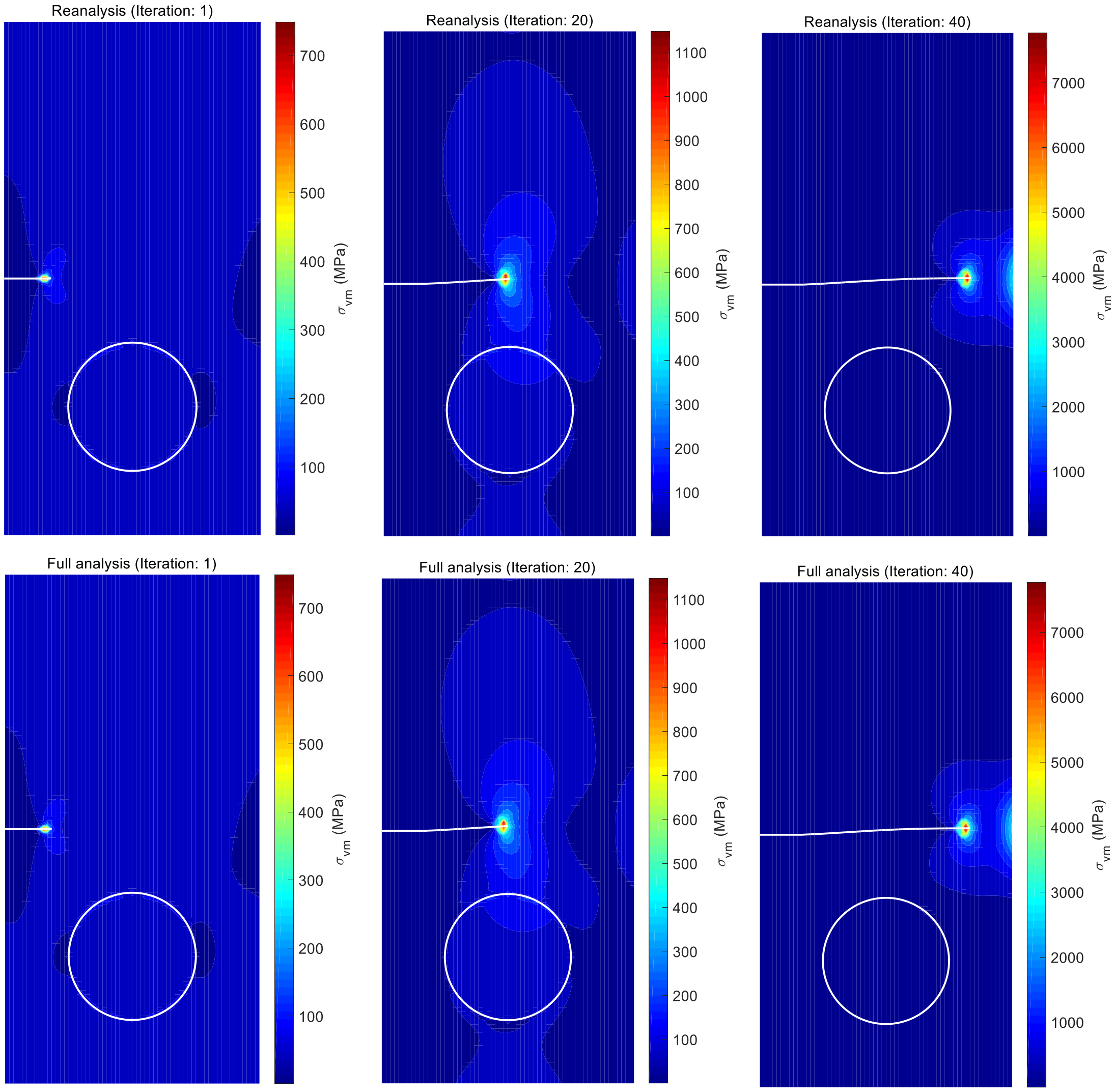

Fig. 15 Von Mises stress contour of edge crack in a plate with a circular inclusion

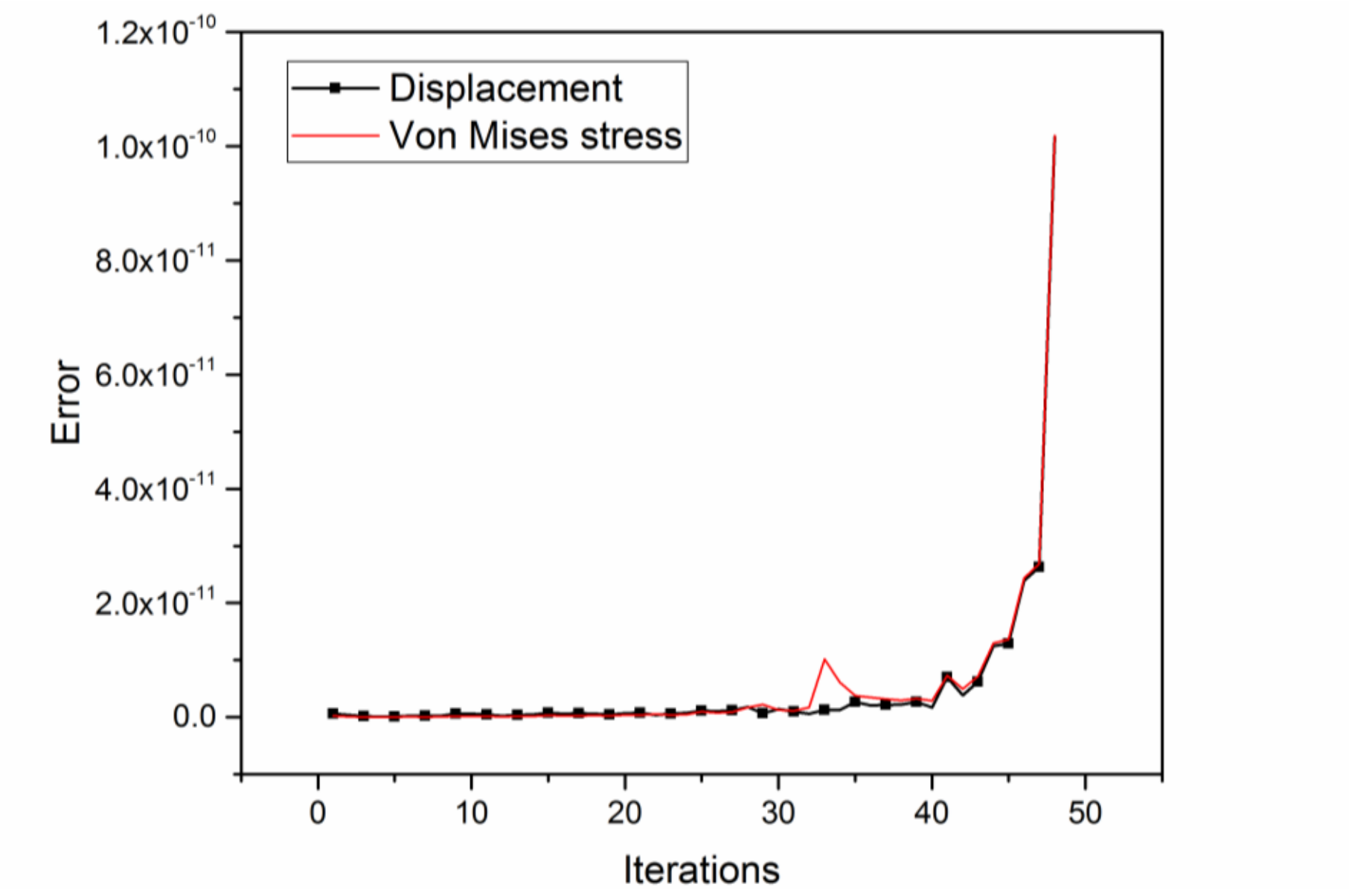

Fig. 16 The errors of each iterative step

Tab. 3 Performance comparison of edge crack in a plate with a circular inclusion

| Computational Time/*s* | | Average Errors | |
|---|---|---|---|
| DUR | Full analysis | Displacement | Von Mises Stress |
| 4.078 | 17.203 | 1.0534e-11 | 1.0951e-11 |

Tab. 4 Crack tip coordinates of edge crack in a plate with a circular inclusion

| Iterative steps | Crack tip coordinates of DUR Method | | Crack tip coordinates of full analysis | |
|---|---|---|---|---|
| | X-coordinate | Y-coordinate | X-coordinate | Y-coordinate |
| 0 | 0.000 | 60.000 | 0.000 | 60.000 |
| | 10.000 | 60.000 | 10.000 | 60.000 |
| 1 | 10.987 | 59.842 | 10.987 | 59.842 |
| 5 | 14.941 | 59.240 | 14.941 | 59.240 |
| 10 | 19.828 | 58.188 | 19.828 | 58.188 |
| 15 | 24.647 | 56.857 | 24.647 | 56.857 |
| 20 | 29.401 | 55.308 | 29.401 | 55.308 |
| 25 | 34.135 | 53.698 | 34.135 | 53.698 |
| 30 | 38.956 | 52.383 | 38.956 | 52.383 |
| 35 | 43.911 | 51.754 | 43.911 | 51.754 |
| 40 | 48.908 | 51.594 | 48.908 | 51.594 |
| 45 | 53.908 | 51.552 | 53.908 | 51.552 |
| 49 | 55.908 | 51.545 | 55.908 | 51.545 |

### 4.3 Center crack in a plate with a circular inclusion and a hole

A center crack in a plate with a circular inclusion and a hole is considered as shown in Fig. 17. The material of plate is aluminum 7075-T6, the Young's modulus $E_1 = 7.17\times10^4 MPa$, the Poisson's ratio $\nu = 0.33$, and a plane strain state is assuming. The material of inclusion is considered as carbon fiber reinforced composite and assume that the Young's modulus $E_2 = 2.1\times10^5 MPa$, the Poisson's is the same as the plate. The initial crack length $a_0 = 10mm$, and the uniformed load $q = 50N/mm$. The increment of propagation $\Delta a = 1mm$. Then, the DUR method is used to solve this problem and only the first iterative step needs to be calculated by the full analysis method while other iterative steps should be predicted efficiently by the DUR method.

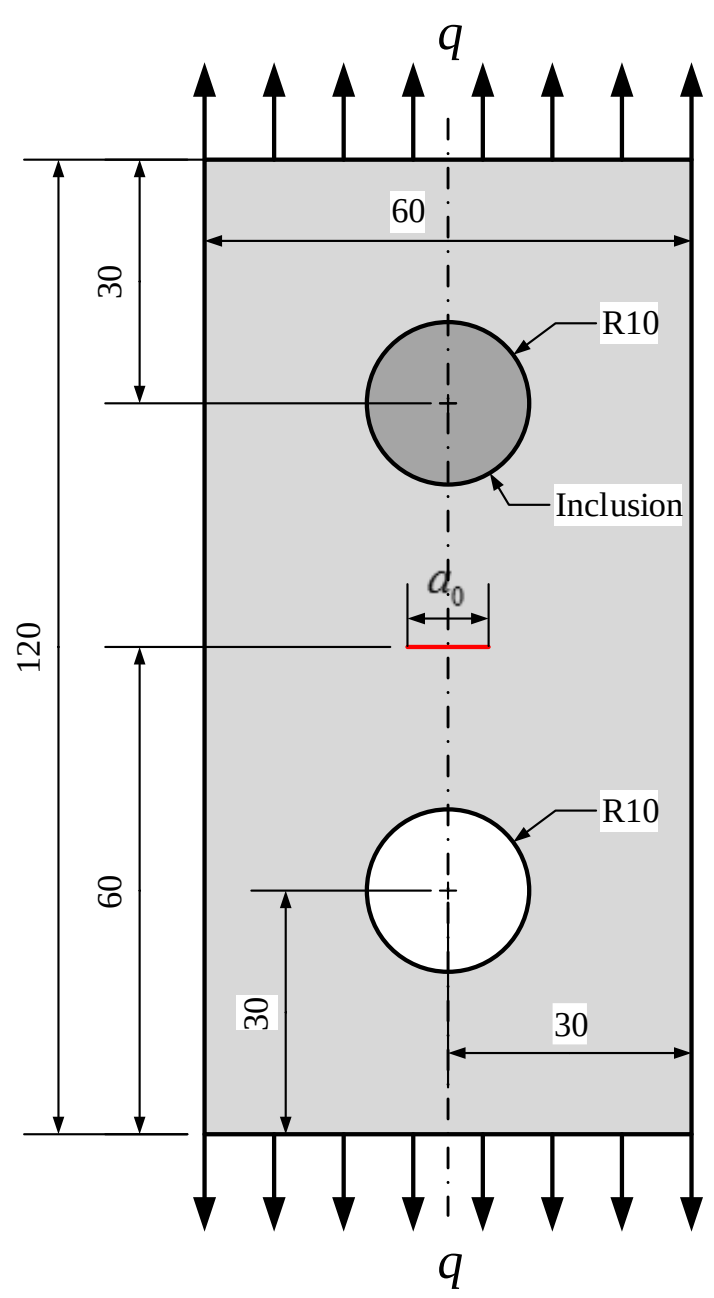

Fig. 17 The geometry of center crack in a plate with a circular inclusion and a hole

The displacement and Von Mises stress contours of iterative step 1, 10 and 24 are shown as Fig. 18 and Fig. 19. It is obviously that the result of the DUR method is very close to the full analysis method. Furthermore, an error list of some selected iterative steps are shown as Fig. 20. It can be found that the DUR method is accurate because the maximum of error is $4.0\times10^{-12}$. Moreover, the computational costs of the DUR and full analysis are also listed in Tab. 5. It shows that the efficiency of the DUR method is much higher than full analysis method. Besides this, the crack tip coordinates are also listed in Tab. 6 and the result of DUR and full analysis are matched.

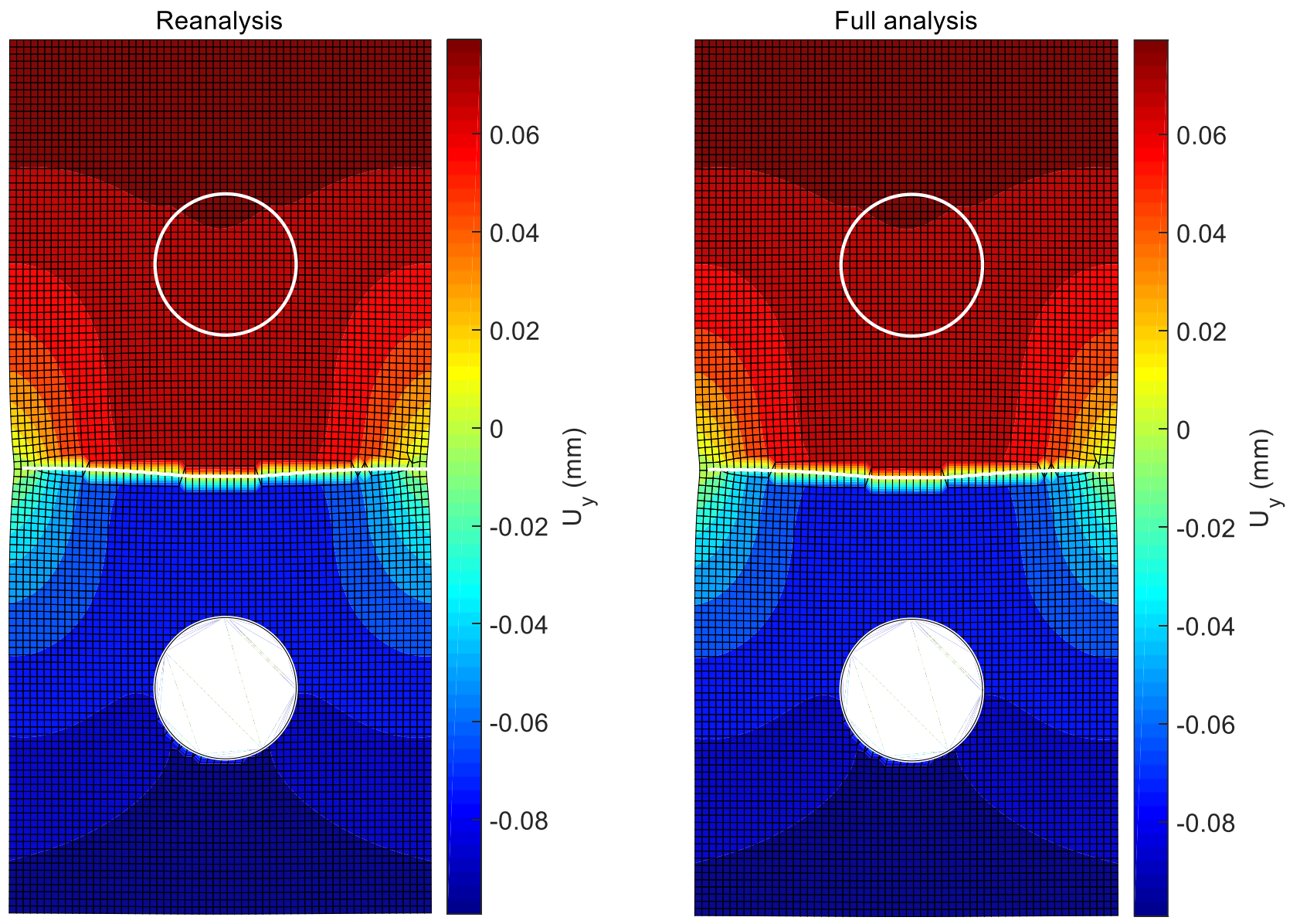

Fig. 18 Displacement contour of center crack in a plate with a circular inclusion and a hole

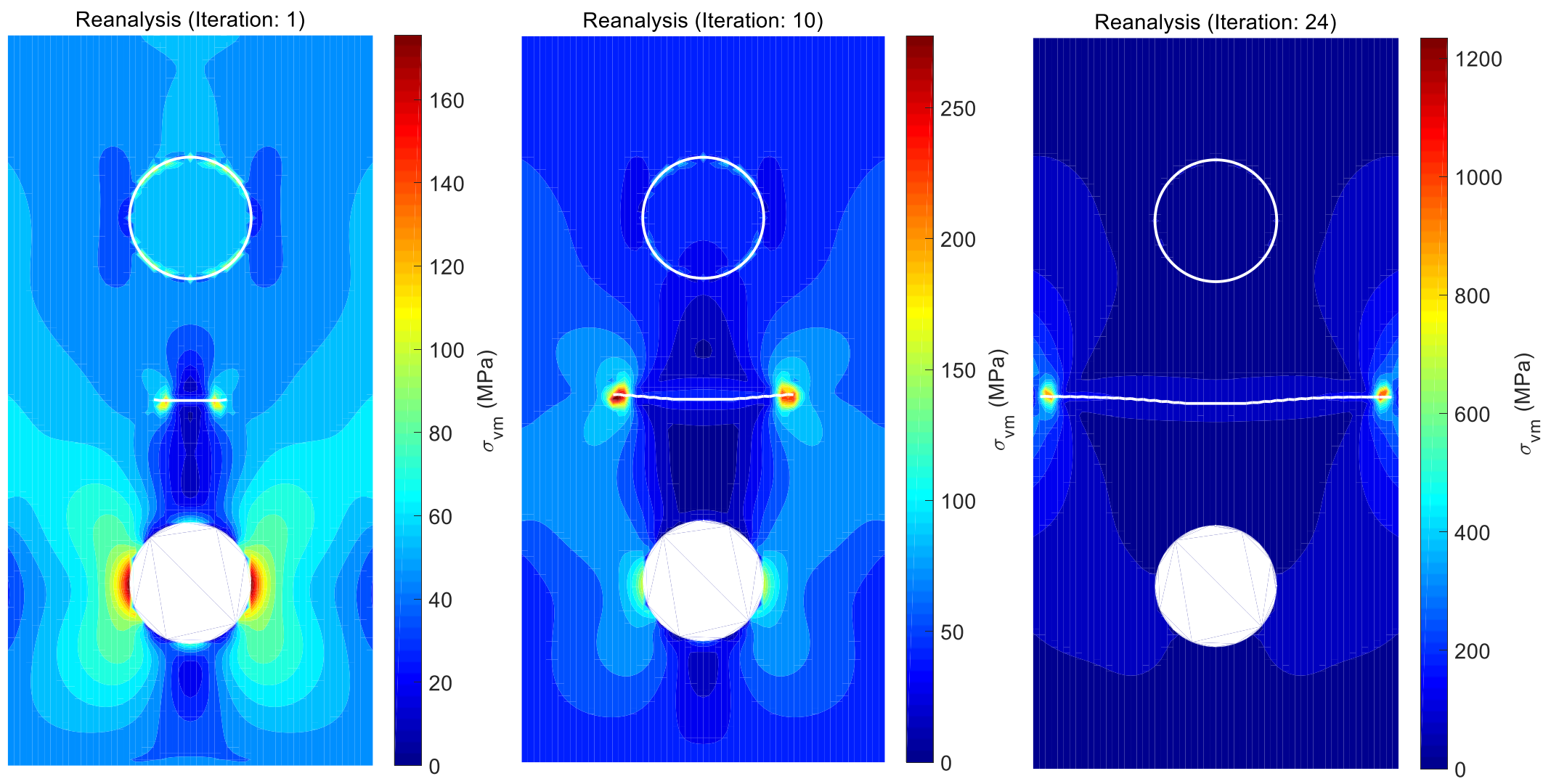

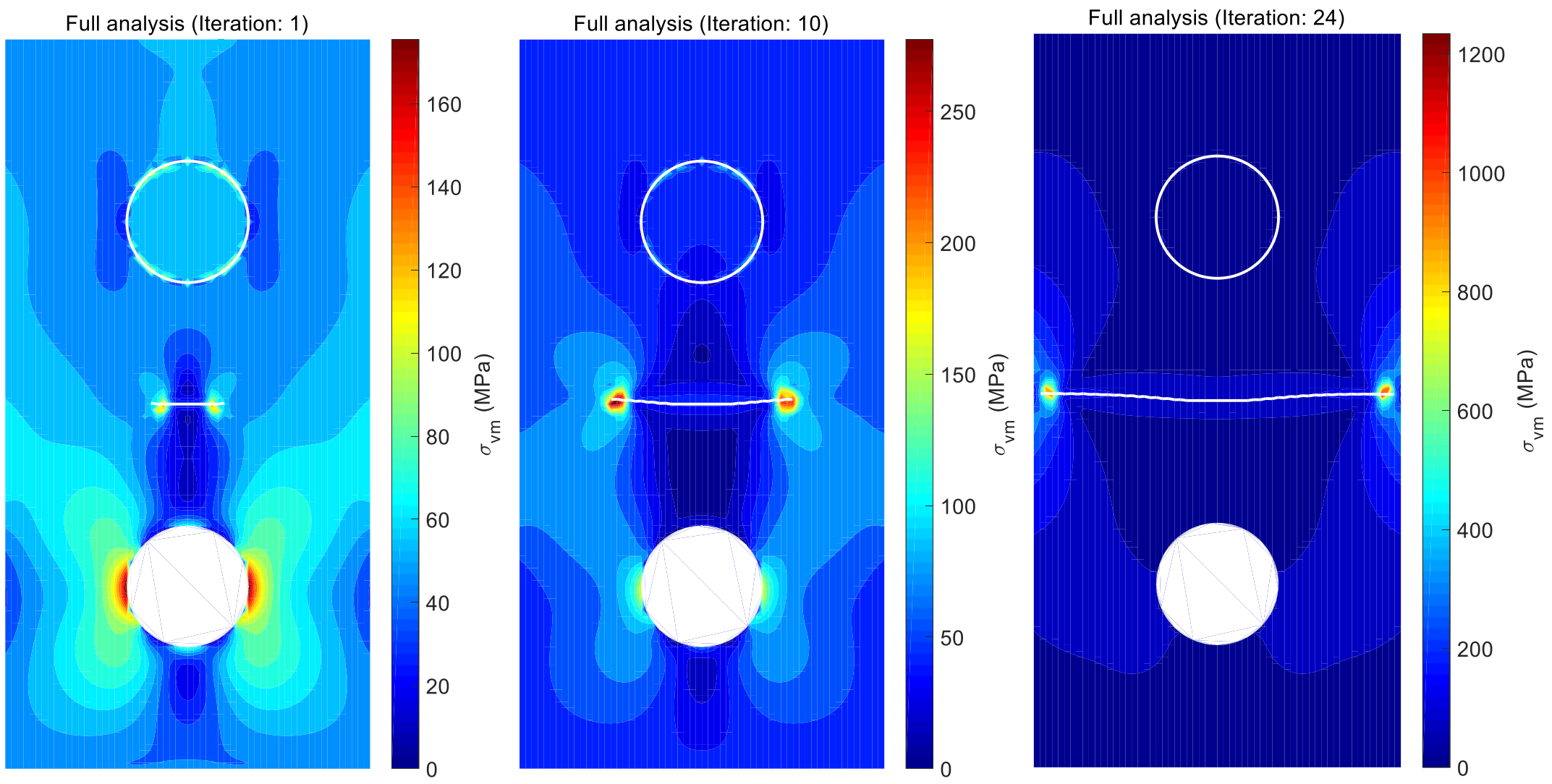

Fig. 19 Von Mises stress contour of center crack in a plate with a circular inclusion and a hole

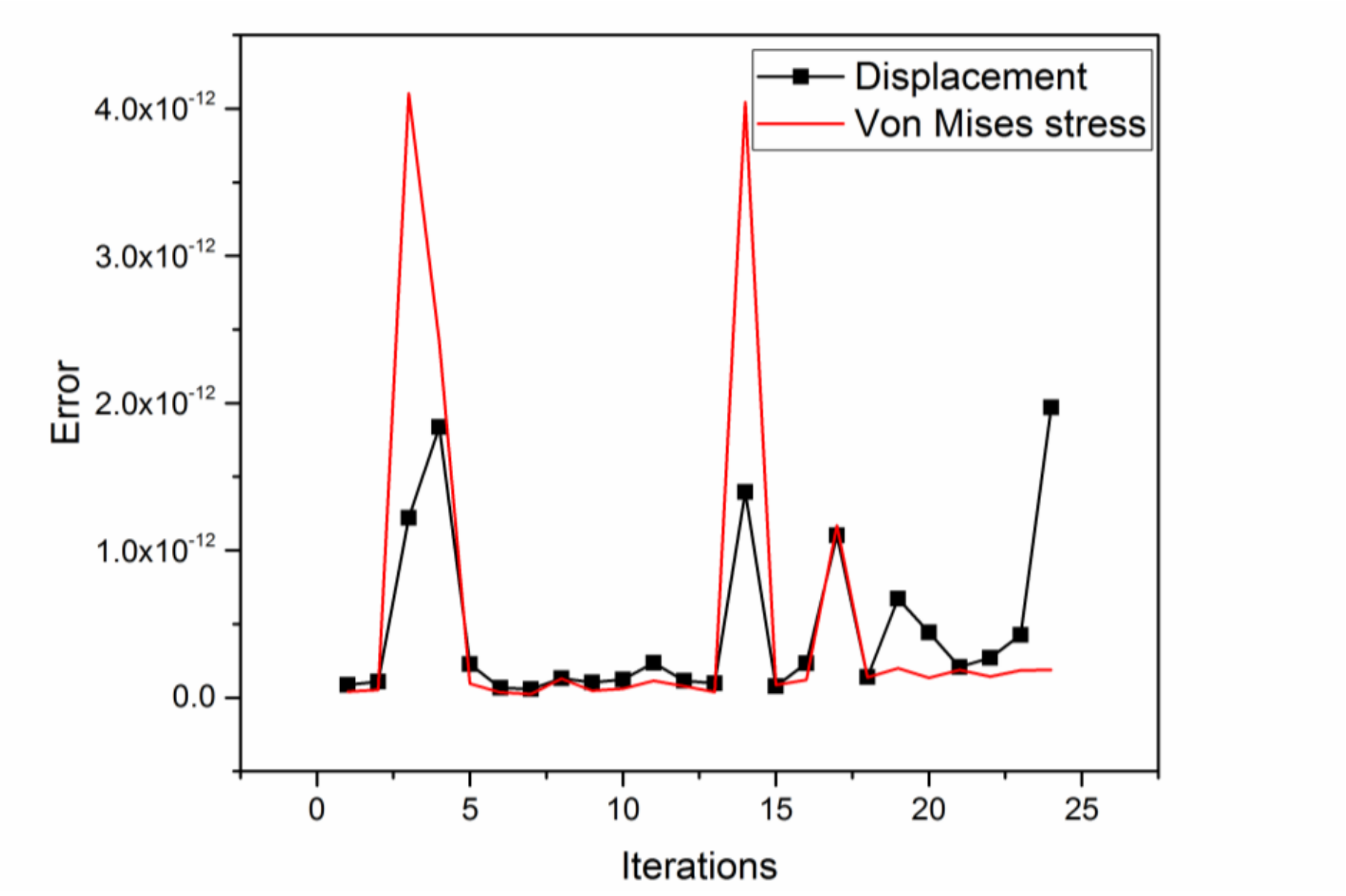

Fig. 20 The errors of each iterative step

Tab. 5 Performance comparison of center crack in a plate with a circular inclusion and a hole

| Computational Time/s | | Average Errors | |
|---|---|---|---|
| DUR | Full analysis | Displacement | Von Mises Stress |
| 2.3438 | 8.4375 | 4.9236e-13 | 5.6032e-13 |

Tab. 6 Crack tip coordinates of center crack in a plate with a circular inclusion and a hole

| Iterative steps | Crack tip coordinates of DUR Method | | Crack tip coordinates of full analysis | |
|---|---|---|---|---|
| | Crack tip 1 | Crack tip 2 | Crack tip 1 | Crack tip 2 |
| 0 | (25.000, 60.000) | (35.000, 60.000) | (25.000, 60.000) | (35.000, 60.000) |
| 1 | (24.012, 60.152) | (35.988, 60.152) | (24.012, 60.152) | (35.988, 60.152) |
| 5 | (20.058, 60.511) | (39.942, 60.510) | (20.058, 60.511) | (39.942, 60.510) |
| 10 | (15.107, 60.748) | (44.894, 60.746) | (15.107, 60.748) | (44.894, 60.746) |
| 15 | (10.133, 61.037) | (49.864, 61.032) | (10.133, 61.037) | (49.864, 61.032) |
| 20 | (5.150,61.077) | (54.858, 61.037) | (5.150,61.077) | (54.858, 61.037) |
| 24 | (1.163, 61.116) | (58.857, 61.019) | (1.163, 61.116) | (58.857, 61.019) |

### 4.4 Accuracy and efficiency comparison

Three numerical examples have been tested in this section and it can be found that the DUR is an accurate and efficient and method. Moreover, a representative case has been calculated by the DUR method under different computational scales from 1000 to 100,000 to fully investigate the efficiency of the DUR method. The log-log plots of comparison results are shown in Fig. 21, Fig. 22 and the error analysis is also shown.

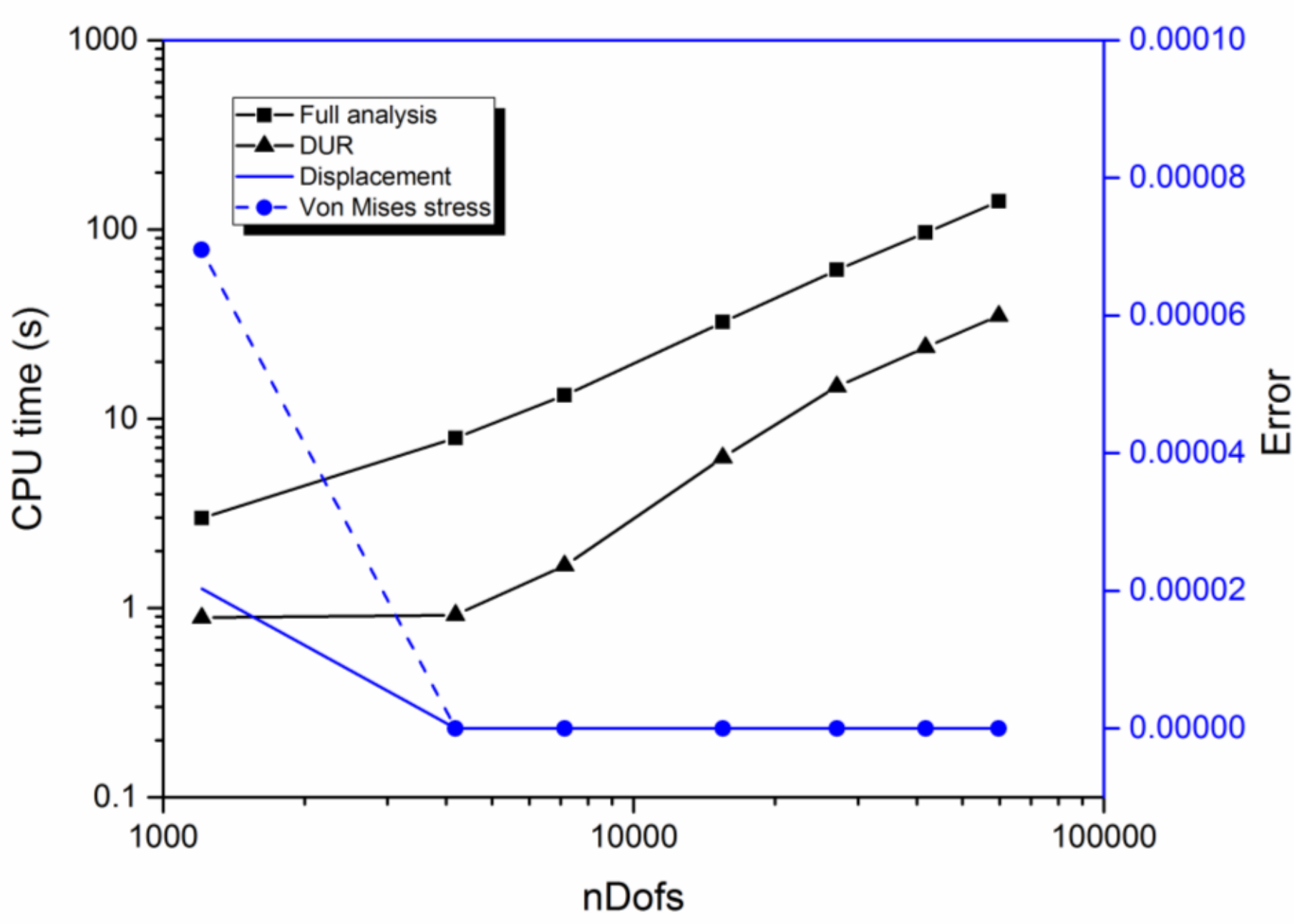

Fig. 21 The comparison of computational cost used in solving equilibrium equations

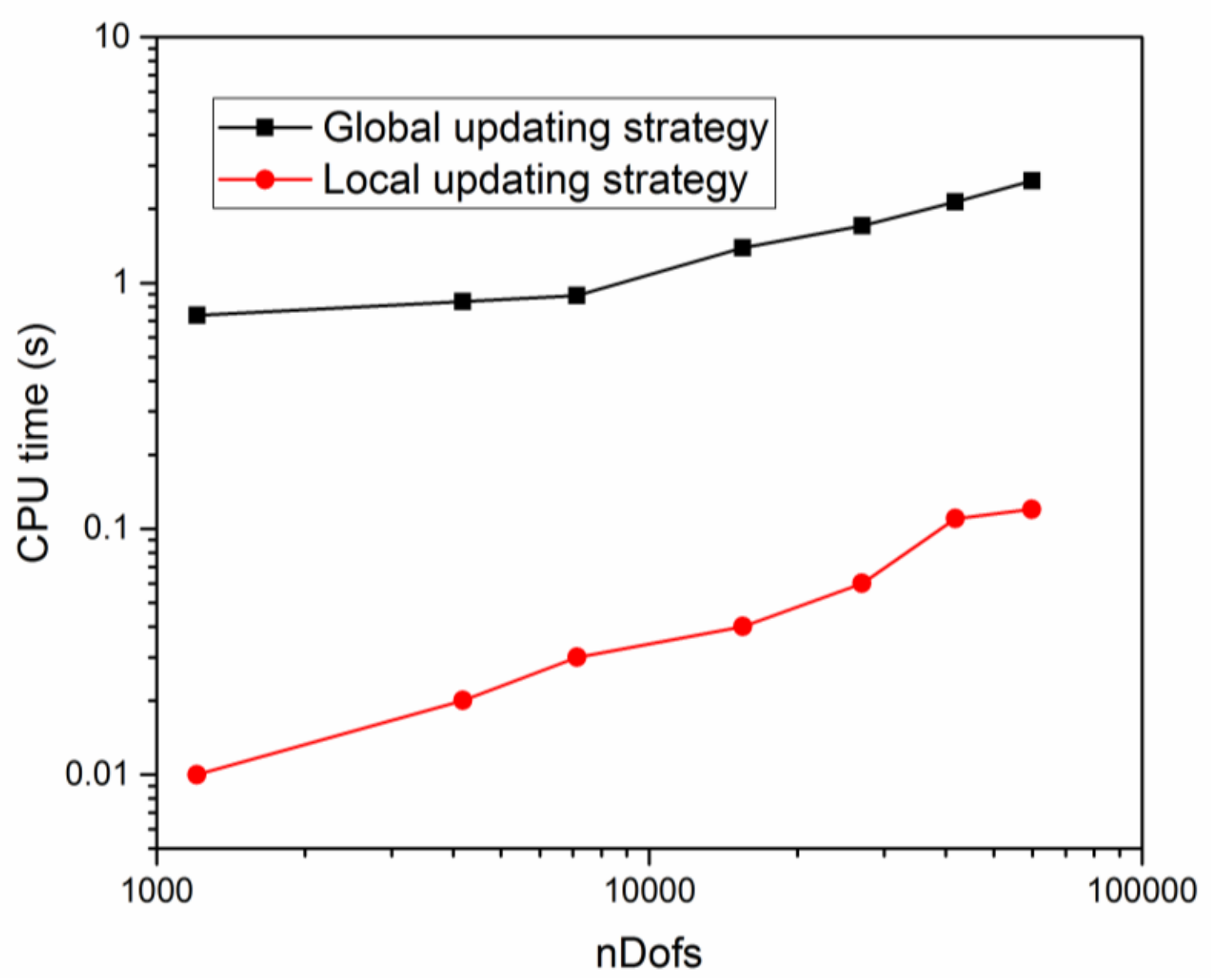

Fig. 22 The comparison of computational cost used in stiffness matrix updating

It can be found that the efficiency of the DUR method is much higher than the full analysis method, and the advantage is more obvious for large scale problems. It is obvious that the accuracy should be improved significantly with the increase of DOFs. Accord to Fig. 21, it can be observed that the accuracy of the DUR method is very high, and the DUR can be regarded as an exact method. Moreover, the computational cost of stiffness matrix updating was also plotted, and Fig. 22 shows that the local updating strategy largely reduces the computational cost of stiffness matrix updating than traditional global updating strategy.

## 5 Conclusions

In this study, the DUR method is proposed for crack propagation. The DUR method consists of three strategies: local stiffness matrix updating, decomposed updating reanalysis, local Cholesky factorization updating strategies. Considering the characteristic of local change of stiffness matrix during X-FEM iterative procedure, the local stiffness matrix updating method can achieve the modified stiffness matrix quickly, and the local Cholesky factorization updating strategy is used to guarantee the accuracy and efficiency of the DUR method. More importantly, the decomposed updating reanalysis strategy is suggested to improve the efficiency of solving the equilibrium equations significantly. Therefore, the DUR method not only reduces the computational cost of solving equilibrium equations but also saves computational cost of stiffness matrix assembling.

Numerical examples show that the DUR method is accurate for the crack propagation. For both the edge and center crack propagation, the accuracy of DUR method is very high. The log-log plots show that the efficiency of DUR method is much more higher than full analysis, and the advantage is more obvious for large scale problems. Moreover, compared with other reanalysis method, the comparisons of stress between DUR and full analysis are made, and the stress can be obtained accurately and efficiently.

**Acknowledgements** This work has been supported by Project of the Key Program of National Natural Science Foundation of China under the Grant Numbers 11572120, National Key Research and Development Program of China 2017YFB0203701.